\documentclass[11pt]{article}

\usepackage{amsmath,amssymb,bm,amsfonts,amsthm,latexsym,url,tabls,multirow}

\theoremstyle{plain}      \newtheorem*{thm*}{Theorem}
                          \newtheorem{thm}{Theorem}[section]
                          \newtheorem{lem}[thm]{Lemma}
                          \newtheorem{prop}[thm]{Proposition}
                          \newtheorem{cor}[thm]{Corollary}     
\theoremstyle{remark}     \newtheorem*{rem}{Remark}
                          \newtheorem*{rems}{Remarks}
\theoremstyle{definition} \newtheorem{example}{Example}
                          \newtheorem*{example*}{Example 1 continued}
                          \newtheorem{fact}{Fact}

    \date{}

\numberwithin{equation}{section}

\usepackage[numbers]{natbib}
\usepackage{varioref}              
\labelformat{section}{Section~#1}  
\labelformat{figure}{Figure~#1}
\labelformat{table}{Table~#1}

\newcommand\x{\mathcal{X}}

\begin{document}

\title{On adding a list of numbers (and other one-dependent
  determinantal processes)}

\author{Alexei Borodin\footnote{\textit{Department of Mathematics, Caltech}}\and
          Persi Diaconis\footnote{\textit{Departments of Mathematics
              and Statistics, Stanford}}\and
          Jason Fulman\footnote{\textit{Department of Mathematics, USC}}}

\date{Version of April 6, 2009}

\maketitle\thispagestyle{empty}

\begin{abstract}
  Adding a column of numbers produces ``carries'' along the way. We show
  that random digits produce a pattern of carries with a neat
  probabilistic description: the carries form a one-dependent
  determinantal point process. This makes it easy to answer natural
  questions: How many carries are typical? Where are they located?
  We show that many further examples, from combinatorics, algebra and group
  theory, have essentially the same neat formulae, and that any one-dependent
  point process on the integers is determinantal. The examples give
  a gentle introduction to the emerging fields of one-dependent and
  determinantal point processes.
\end{abstract}

\section{Introduction}\label{sec1}

Consider the task of adding a single column of digits:

\begin{equation*}\begin{array}{rrr}
7&\cdot&7\\
9&\cdot&6\\
4&&0\\
8&\cdot&8\\
3&&1\\
6&&7\\
1&\cdot&8\\
6&\cdot&4\\
8&&2\\ \cline{1-1}
52&&
\end{array}\end{equation*}
We have put a dot to the right of a digit whenever the succeeding addition
leads to a carry. The remainder (mod 10) is written to the right of the
dot. Thus $7+9=16$ results in a dot by the 7 and a remainder of 6. At the
end, the total is found by adding the number of dots (here, five) and
appending the final remainder (here, 2) for a total of 52. The dots are a
standard bookkeeping device used (for example) in the Trachtenberg system
of speed addition \citep{cm}.

How many carries are typical and how are they distributed? Common
sense suggests that about half the places have carries and that if a
carry occurs, it is less likely that the next addition gives a
carry. Investigating these questions when digits are chosen uniformly
at random leads to interesting mathematics.

As additional motivation, look at the remainder column in the example
above. Observe that there is a carry on the left (and so a dot) if and
only if there is a descent on the right (a sequence $x_1,x_2,\dots$
has a descent at $i$ if $x_i>x_{i+1}$). Further, if the original
column of digits is independent and uniformly distributed, so is the
remainder column. Thus, the results about carries give the
distribution of the descent pattern in a random sequence of digits.

\ref{sec2} derives the distribution theory of carries by elementary
arguments. There are nice formulae for the chance of any pattern of
carries, and the process is one-dependent so a variety of probabilistic
limit theorems are available.

\ref{sec3} introduces determinantal point processes, shows that the
carries process is determinantal, and illustrates how standard tools for
determinantal point processes apply to carries.

\ref{newsec4} reviews the literature on stationary one-dependent
processes, shows that all of these are determinantal and further that
essentially all the neat formulae for carries and descents have versions
for any stationary one-dependent process. Connections with symmetric
function theory are developed. A large class of examples arising from the
work of Polishchuk and Positselski \cite{pp} on Koszul algebras is shown
to yield many natural examples (including the original carries process).

\ref{sect5} gives combinatorial examples: descents in permutations from
both uniform and non-uniform distributions, the connectivity set of a
permutation, and binomial posets.

\ref{sect6} generalizes from adding numbers to multiplying random elements
in a finite group. For central extensions (e.g., the quaternions), there
is again a carries process that is stationary, one-dependent, and
determinantal. Thus explicit formulae and limit theorems are available.

\ref{sect7} contains proofs of some of our more technical results. It
shows that any one-dependent (possibly non-stationary) point process on
the integers is determinantal. It also constructs a family of
one-dependent determinantal processes (generalizing many examples in
earlier sections), and computes its correlation kernel.

\section{Probability theory for carries and descents}\label{sec2}

Throughout this section we work base $b$ and so with the alphabet
$\mathcal{B}=\{0,1,\dots,b-1\}$. Let $B_1,B_2,\dots,B_n$ be a sequence of
randomly chosen elements of $\mathcal{B}$ (independent and identically
distributed). There is a \textit{descent at $i$} if $B_i>B_{i+1},\
1\leq i\leq n-1$. Let $X_i$ be $1$ or $0$ as there is a descent at $i$
and $D=\{i:X_i=1\}$ be the descent set of $B_1,B_2,\dots,B_n$.

The following probabilistic facts are elementary. They are stated for
carries but verification is easier using descents.
\begin{fact}
(Single carries) \quad For any $i\in[n-1]$,
\begin{equation*}
P(X_i=1)=\dfrac12-\dfrac1{2b}=\frac{\binom{b}{2}}{b^2}.
\end{equation*}
Thus, when $b=10$, the chance of a carry is $.45$. When $b=2$, the chance
of a carry is $.25$. For any base, Var$(X_i)=\tfrac14-\tfrac1{4b^2}$.
\end{fact}
\begin{fact}
(Runs of carries) \quad For any $i$ and $j$ with $1\leq
  i<i+j\leq n$,
\begin{equation*}
P\left(X_i=X_{i+1}=\dots=X_{i+j-1}=1\right)=\binom{b}{j+1}\bigg/b^{j+1}.
\end{equation*}
Thus, a run of $b$ or more carries is impossible. Further,
\begin{equation*}
\text{Cov}\left(X_i,X_{i+1}\right)=E\left(X_iX_{i+1}\right)-E(X_i)E(X_{i+1})=-\dfrac1{12}\left(1-\dfrac1{b^2}\right).
\end{equation*}
\end{fact}
\begin{fact}
(Stationary one-dependence) \quad The distribution of
  $\{X_i\}_{i\in[n-1]}$ is stationary: for $J\subseteq[n-1],\ i\in[n-1]$
  with $J+i\subseteq[n-1]$, the distribution of $\{X_j\}_{j\in J}$ is the
  same as the distribution of $\{X_j\}_{j\in J+i}$. Further, the
  distribution of $\{X_i\}_{i\in[n-1]}$ is one-dependent: if
  $J\subseteq[n-1]$ has $j_1,j_2\in J\Rightarrow|j_i-j_2|>1$, then
  $\{X_j\}_{j\in J}$ are jointly independent binary random variables
  with $P(X_j=1)=\tfrac12-\tfrac1{2b}$. The literature on
  $m$-dependent random variables is extensive. In particular, a
  classical central limit theorem \cite{hr} shows the following:

\begin{thm} For $n \geq 2$, the total
number of carries $T_{n-1}=X_1+\dots+X_{n-1}$ has mean
$(n-1)(\tfrac12-\tfrac1{2b})$, variance $\frac{n+1}{12}(1-\tfrac1{b^2})$
and, normalized by its mean and variance, $T_{n-1}$ has a standard normal
limiting distribution for $n$ large. \end{thm}
 {\it Remarks.}
\begin{enumerate} \item An $O(n^{-1/2})$ error bound in this central limit
theorem can be proved using the dependency graph approach to normal
approximation by Stein's method \cite{cs}, \cite{der}.

\item
Stationary pairwise independent processes can fail to obey the central
limit theorem \cite{brad}, \cite{jan}.
\end{enumerate}
\end{fact}
\begin{fact}
($k$-point correlations) \quad For $A\subseteq[n-1]$, let
\[\rho(A)=P\{X_i=1\text{
  for }i\in A\}.\] For $|A|=k$, the $\rho(A)$ are called  $k$-point
correlations for the point process $X_1,\dots,X_{n-1}$. They are basic
descriptive units for general point processes \cite{dv}. For the carries
and descent process, they are simple to describe. Break $A\subseteq[n-1]$
into disjoint, non-empty blocks of consecutive integers $A=A_1\cup
A_2\cup\dots\cup A_k$. Thus
$A=\{2,3,5,6,7,11\}=\{2,3\}\cup\{5,6,7\}\cup\{11\}$. From Facts 2 and 3
above,
\begin{equation*}
\rho(A)=\prod_{i=1}^k \left[ \binom{b}{a_i+1}\bigg/b^{a_i+1} \right]
\qquad\text{if }A=\bigcup_{i=1}^kA_i\text{ with }|A_i|=a_i.
\end{equation*}
\end{fact}
\begin{fact}
(Determinant formula) \quad Let
$\epsilon_1,\epsilon_2,\dots,\epsilon_{n-1}$ be a fixed sequence in
$\{0,1\}$. Then
\begin{equation*}
P\left\{X_1=\epsilon_1,\dots,X_{n-1}=\epsilon_{n-1}\right\}=\dfrac1{b^n}
\cdot \text{det}\binom{s_{j+1}-s_i+b-1}{b-1}.
\end{equation*}
Here, if there are exactly $k$ $1$'s in the $\epsilon$-sequence at
positions $s_1<s_2<\dots<s_k$, the determinant is of a $(k+1)\times(k+1)$
matrix with $(i,j)$ entry $\binom{s_{j+1}-s_i+b-1}{b-1}$ for $0\leq
i,j\leq k$ with $s_0=0,\ s_{k+1}=n$.
\end{fact}
\begin{example}
If $n=8,\ \epsilon_1=1,\ \epsilon_2=\epsilon_3=\epsilon_4=0,\ \epsilon_5=1,\ \epsilon_6=\epsilon_7=0$, the matrix is
\begin{equation*}\begin{pmatrix}
\binom{1+b-1}{b-1}&\binom{5+b-1}{b-1}&\binom{8+b-1}{b-1}\\
1&\binom{4+b-1}{b-1}&\binom{7+b-1}{b-1}\\
0&1&\binom{3+b-1}{b-1}
\end{pmatrix}.
\end{equation*}
When $b=2$, this is
$\left(\begin{smallmatrix}2&6&9\\1&5&8\\0&1&4\end{smallmatrix}\right)$
with determinant $9$. Thus the chance of a carry at exactly positions 1
and 5 when adding eight binary numbers is $9/2^8\doteq.03516$. When
$b=10$, this chance becomes $1,042,470/10^8\doteq.0104$. This is not so
much smaller than in the binary case. \label{ex21}
\end{example}
\begin{proof} (Of fact 5)
We follow Stanley \cite[p.~69]{sta2}. Let $\alpha_n(S)$ be the number of
sequences with descent set contained in $S$ and $\beta_n(S)$ the number of
sequences with descent set equal to $S$. Then
\begin{equation*}
\alpha_n(S)=\sum_{T\subseteq S}\beta_n(T),\qquad
\beta_n(S)=\sum_{T\subseteq S}(-1)^{|S-T|}\alpha_n(T).
\end{equation*}
There is a simple formula for $\alpha_n(S)$: the number of weakly
increasing sequences with entries in $\{0,1,\dots,b-1\}$ of length $l$
is $\binom{l+b-1}{b-1}$ using the usual stars and bars argument (place
$b-1$ bars into $l+b-1$ places, put symbol $0$ to the left of the first
bar, symbol $1$ between the first and second bar, and so on, with symbol
$b-1$ to the right of the last bar). Then
\begin{equation*}
\alpha_n(S)=\binom{s_1+b-1}{b-1}\binom{s_2-s_1+b-1}{b-1}\dots\binom{n-s_k+b-1}{b-1},
\end{equation*}
because any such sequence of length $n$ is constructed by
concatenating nondecreasing sequences of length $s_1,
s_2-s_1,\dots,n-s_k$. Therefore
\begin{equation*}
\beta_n(S)=\sum_{1 \leq i_1<i_2<\dots<i_j\leq k} (-1)^{k-j}
f(0,i_1)f(i_1,i_2)\dots f(i_j,k+1)
\end{equation*}
with $f(i,j)=\binom{s_j-s_i+b-1}{b-1}$ (with $s_0=0,\ s_{k+1}=n$). From
Stanley's discussion, $\beta_n(S)$ is the determinant of the
$(k+1)\times(k+1)$ matrix with $(i,j)$ entry $f(i,j+1),\ 0\leq i,j\leq k$.
\end{proof}
\begin{example}
From the above development, the chance that the sum of $n$ base $b$ digits
has \textit{no} carries is $\binom{n+b-1}{b-1}/b^n$. When $n=8,\ b=2$,
this is $9/2^8\doteq.03516$. When $n=8,\ b=10$, this is
$\binom{17}{9}/10^8\doteq.00024$. We feel lucky when adding eight numbers
with no carries. These calculations show that we are lucky indeed.
\label{ex22}
\end{example}

\section{Determinantal point processes}\label{sec3}

Let $\x$ be a finite set. A point process on $\x$ is a probability
measure $P$ on the $2^{|\x|}$ subsets of $\x$. For example, if
$\x=\{1,2,\dots,n-1\}$, a point process is specified by recording
where the carries occur in adding $n$ base $b$ numbers as in
\ref{sec2}. One simple way to specify $P$ is via its so-called
correlation functions. For $A\subseteq\x$, let
\begin{equation*}
\rho(A)=P\{S:S\supseteq A\}.
\end{equation*}
The collection of numbers $\{\rho(A)\}$ uniquely determine $P$ using
inclusion/exclusion. A point process is \textit{determinantal} with
kernel $K(x,y)$ if
\begin{equation*}
\rho(A)=\text{det}\left(K(x,y)\right)_{x,y\in A}.
\end{equation*}
On the right is the determinant of the $|A|\times|A|$ matrix with
$(x,y)$ entry $K(x,y)$ for $x,y\in A$.

Determinantal point processes were introduced by Macchi \cite{ma} to model
the distribution of fermions. See \cite{dv} for a textbook account,
\cite{sos} for a survey, \cite{hkpv} for probabilistic developments and
\cite{ly} for many combinatorial examples. We warn the reader that these
last references essentially deal with symmetric kernels $K(x,y)$ whereas
all of our examples involve non-symmetric kernels.

There has been an explosive development of determinantal point
processes over the past few years. Some reasons for this include:
\begin{enumerate}
\item[(a)] A raft of natural processes (including all those in the present
  paper) turn out to be determinantal. These include: the eigenvalues
  of various ensembles of random matrices (Dyson \cite{dy1,dy2,dy3},
  Mehta \cite{meh}); the presence or absence of edges if a random
  spanning tree is chosen in a graph (Burton--Pemantle \cite{bp},
  Lyons \cite{ly}); random tilings and growth models (Johansson \cite{j2});
  the structure of random partitions chosen from
  a variety of measures (e.g., the Plancherel measure on the symmetric
  group $S_n$) (Borodin--Okounkov--Olshanski
  \cite{boo}); random dimers on bipartite planar graphs (Kenyon \cite{ke});
  and the zeros of a random analytic function
  $\sum_{n=0}^\infty a_n z^n$ with i.i.d. complex Gaussian coefficients
  $a_n$ \cite{pv}. The list goes on extensively.

\item[(b)] Specifying a kernel may be much easier than specifying a measure
  on all $2^{|\x|}$ subsets. Further, the correlation functions allow
  easy computation of quantities of interest. For example, if $X_x$ is
  $1$ or $0$ as the random set includes $x$ or not:
\begin{align*}
E(X_x)&=K(x,x),\\
\text{Cov}(X_x,X_y)&=\text{det}\begin{pmatrix}
K(x,x)&K(x,y)\\K(y,x)&K(y,y)
\end{pmatrix}-K(x,x)K(y,y)\\
&=-K(x,y)K(y,x).
\end{align*}
If $K(x,y)=K(y,x)$, then the correlation is $\leq 0$, a distinctive
feature of determinantal point processes.
\item[(c)] If the matrix $(K(x,y))_{x,y\in\x}$ has all real eigenvalues
  $\{\lambda_x\}_{x\in\x}$, many theorems become available for $N$, the
  total number of points in a realization of the process. These
  include:
\begin{itemize}
\item $N$ is distributed as a sum $\sum_{x\in\x}Y_x$ with $\{Y_x\}$
  independent $0/1$ random variables having $P(Y_x=1)=\lambda_x$.
\item Let $\mu$ and $\sigma$ denote the mean and variance of
  $N$. These are available in terms of $K$ from (b) above. Then the
  following refined form of the central limit theorem holds:
\begin{equation*}
\left|P\left\{\frac{N-\mu}{\sigma}\leq
    x\right\}-\dfrac1{\sqrt{2\pi}}\int_{-\infty}^x e^{-t^2/2}dt\right|\leq
    \frac{.80}{\sigma}
\end{equation*}
\item The probability density $P(n)=P(N=n)$ is a Polya-frequency
  function. In particular the density is unimodal. Further,
  $|\text{mode}-\text{mean}|\leq1$.
\item There is a useful algorithm for simulating from $K$ (here
  assuming $K(x,y)=K(y,x)$).
\end{itemize}
Most of these properties follow from the fact that the generating
function of $N$,
\begin{equation*}
E(x^N)=\text{det}\left(I+(x-1) K\right)
\end{equation*}
has all real zeros. See Pitman \cite{Pi} or Hough, et al. \cite{hkpv} for
details.

\item[(d)] If $K_n$ is a sequence of such kernels converging to a kernel
  $K$, then, under mild restrictions, the associated point processes
  converge. See Soshnikov \cite{sos} or Johansson and Nordenstam \cite{jn} for
  precise statements and remarkable examples. A simple example is
  given in \ref{sect5} (Example \ref{exam41}).
\end{enumerate}
A main result of our paper is that the carries and descent processes, and
indeed any one-dependent point processes on the integers, is
determinantal. The special case of carries is stated here; a proof of the
general case is in \ref{sect7}.
\begin{thm}
The point process of carries and descents $P_n$ given in \ref{sec2} is
determinantal with correlation kernel $K(i,j)=k(j-i)$ where
\begin{equation*}
\sum_{m\in \mathbb{Z}}k(m)t^m=\dfrac1{1-(1-t)^b}.
\end{equation*}
\label{thm31}
\end{thm}
\begin{example}
Take $b=2$. Then $\tfrac1{1-(1-t)^2}=\tfrac1{2t(1-t/2)}$. Replacing $t$ by
$ct$ changes $k(m)$ to $k(m)c^m$ and $K(i,j)=k(j-i)$ to
$c^{j-i}k(j-i)=c^{-i}K(i,j)c^j$. Conjugating $K$ by a diagonal matrix does
not change determinants or the correlation functions. Thus setting
$\tfrac{t}2=z$, the generating function becomes
\begin{equation*}
\dfrac14\left(\dfrac1{z}+1+z+z^2+\dots\right).
\end{equation*}
When e.g., $n-1=5$, the kernel is
\begin{equation}
\left(\dfrac14\right)\begin{bmatrix}
1&1&1&1&1\\
1&1&1&1&1\\
0&1&1&1&1\\
0&0&1&1&1\\
0&0&0&1&1\\
\end{bmatrix}.
\label{totpos}
\end{equation}
From this, $\rho(A)={\footnotesize\begin{cases}0&\text{if $A$ has two consecutive
    entries}\\ \left(\tfrac14\right)^{|A|}&\text{otherwise}\end{cases}}$.
This is a manifestation of the one-dependence together with the fact
that (for $b=2$), two consecutive carries are impossible. It follows
that the \textit{possible} configurations are all binary strings of
length $n-1$ with no two consecutive $1$'s. This is a standard coding
of the Fibonacci number $F_{n+1}$ (if $1,1,2,3,5$ are the first five
Fibonacci numbers). For example, when $n=4$, the possible
configurations are $000,001,010,100,101$.

The eigenvalues of any finite matrix of the form (\ref{totpos}) are real.
This follows from the fact that these matrices are totally positive, see
\cite{ed}, and the well known fact that the eigenvalues of a totally
positive matrix are real (and nonnegative), see e.g. Corollary 6.6 of
\cite{an}.

The classification of \cite{ed} also implies that the correlation kernels
for $b\ge 3$ are not totally positive (for large enough matrix size), and
we do not know if their eigenvalues are real.

\label{ex31}
\end{example}
\begin{example}
Take $b=3$. A straightforward expansion leads to
\begin{eqnarray*}
\sum_{m=-\infty}^\infty k(m)t^m & = & \dfrac1{1-(1-t)^3}=
\frac{1}{3t \left( 1-t+ \frac{t^2}{3} \right)} \\
& = & \frac{ \left(1+t+\frac{2t^2}{3}+\frac{t^3}{3}+\frac{t^4}{9} \right)}
{3t \left(1+\frac{t^6}{27} \right)}. \label{31}
\end{eqnarray*} Thus
\begin{itemize}
\item $k(n)=0$ for $n<-1$,
\item $k(-1)=\frac{1}{3}$,
\item $k(n)=0$ for $n=4+6j,\ 0\leq j<\infty$. For other $n\geq2$,
\item
  $k(n)=(-1)^{\lfloor\tfrac{n+1}6\rfloor}\left(\tfrac13\right)^{\lfloor \tfrac{n+3}2
  \rfloor} 2^{\delta(n)}$, where $\delta(n)=1$ if $n=1 \mod 6$ and $0$ else.
\end{itemize} For example, the first few values are:
\begin{equation*}\begin{small}\begin{array}{r|rrrrrrrrrrrrrrrrrr}
m&-1&0&1&2&3&4&5&6&7&8&9&10&11&12&13&14&15&16\\\hline
k(m)&\tfrac13&\tfrac13&\tfrac29&\tfrac19&\tfrac1{27}&0&-\tfrac1{3^4}&-\tfrac1{3^4}
&-\tfrac2{3^5}&-\tfrac1{3^5}&-\tfrac1{3^6}&0&
\tfrac1{3^7}&\tfrac1{3^7}&\tfrac2{3^8}&\tfrac1{3^8}&\tfrac1{3^9}&0
\end{array}\end{small}
\end{equation*}

It is amusing to see that the two-variable kernel $K(x,y)=k(y-x)$
reproduces the correlation functions found for this problem in \ref{sec2}
(Fact 4). Thus for $n=7$,
\begin{equation*}
(K(x,y))_{i,j=1}^6= \begin{bmatrix}
1/3 & 2/9 & 1/9 & 1/27& 0 &-1/81\\
1/3& 1/3 & 2/9& 1/9 & 1/27& 0\\
0& 1/3& 1/3 & 2/9 & 1/9& 1/27\\
0&0& 1/3& 1/3& 2/9& 1/9 \\
0&0&0&1/3&1/3&2/9\\
0&0&0&0&1/3&1/3\end{bmatrix}
\end{equation*}
\label{ex32}
\end{example}

\section{One-dependent processes}\label{newsec4}

Let $\{X_i\}_{i=0}^N$ be a binary stochastic process. Here
$0<N\leq\infty$. The process is \textit{one-dependent} if
$\{X_j\}_{j=0}^{i-1},\ \{X_j\}_{j=i+1}^N$ are independent for all
$i$. A natural example: let $\{Y_i\}_{i=0}^{N+1}$ be independent
  uniform random variables on $[0,1]$. Let $h:[0,1]^2\to\{0,1\}$ be a
  measurable function. Then $X_i=h(Y_i,Y_{i+1}),\ 0\leq i\leq N$, is a
  one-dependent process called a \textit{two-block factor}. Much of
  the literature on one-dependence is for stationary
  processes. However, if $Y_i$ are independent but not identically
  distributed or $h=h_i$ depends on $i$, the associated $X_i$ is still
  one-dependent. Our main results in \ref{sect7} are for general
  one-dependent processes. We first describe some results for
  the stationary case.

\subsection*{Stationary one-dependent processes}

A book-length review of this subject by de Valk \cite{vdv} collects
together many of the results, so we will be brief, focusing on later
research. One of the first problems studied is the question of
whether all one-dependent processes are two-block factors.
Counterexamples were found \cite{vdv} though it appears that most
``natural'' one-dependent processes are two-block factors. For
example, let $f\in L^2(\pi)$ be a function on the unit circle. Let
$K(i,j)=\hat{f}(i-j),\ i,j\in\mathbb{Z}$, $\hat{f}$ the Fourier
transform. Macchi \cite{ma} and Soshnikov \cite{sos} show that
$K(i,j)$ is the kernel for a determinantal point process on
$\mathbb{Z}$. Lyons and Steif \cite{ls} and Shirai and Takahashi
\cite{st1,st2} develop remarkable properties of the associated point
process. If $f(\theta)=ae^{-i\theta}+b+ce^{i\theta}$, the point
process is one-dependent and \cite{brom} shows it is a two-block
factor. See also \cite{cf}.

A stationary one-dependent process is determined by the numbers
$a_i=P(X_1=X_2=\dots=X_{i-1}=1),\ a_1=1$. Indeed, if
$f=f_1,\dots,f_{i-1};e=e_{i+1},\dots,e_{n-1}$ are binary strings,
$P(f,0,e)+P(f,1,e)=P(f)P(e)$. From this $P(f,0,e)=P(f)P(e)-P(f,1,e)$. By
similar reductions, the probability of any pattern of occurrences can be
reduced to a polynomial in the $a_i$. For example,
$P(0,0,0)=1-3a_2+a_2^2+2a_3-a_4$. There is a remarkable, simple formula
for this polynomial as a minor of a Toeplitz matrix. We learned this from
\citep[Chap.~7]{pp}. A generalization to the non-stationary case is in
\ref{sect7}.
\begin{thm} \label{newthm41}
For a stationary one-dependent process, let
$a_i=P(X_1=X_2=\dots=X_{i-1}=1),\ a_1=1$. Let $t_1,\dots,t_{n-1}$ be a
binary string with $j$ zeros at positions $S = \{s_1<\cdots<s_k\}
\subseteq[n-1]$. Then \[ P(t_1,\cdots,t_{n-1})=
\det\left(a_{s_{j+1}-s_{i}} \right)_{i,j=0}^{k} .\] Here the determinant
is of a $k+1$ by $k+1$ matrix, and one sets $s_0=0,s_{k+1}=n, a_0=1,$ and
$a_i=0$ for $i<0$.
\end{thm}
\begin{proof}
  Theorem \ref{newthm41} follows from Proposition 2.2 and
  Theorem 4.2 in \citep[Chap.~7]{pp}, after elementary manipulations.
\end{proof}

{\it Remark:} It follows from the discussion on \citep[p.~140]{pp} that if
all of the determinants (for all $n$ and subsets $S$) in Theorem
\ref{newthm41} are non-negative, then a stationary, one-dependent process
with these $a_i$ values exists.

\begin{example} For any $n$, the sequence with $n-1$ ones
has $k=0$ zeros. The relevant matrix is then $1$ by $1$ with entry $a_n$,
giving this probability. \label{newex42}
\end{example}
\begin{example}
To compute $P(0,0,0)$ one applies the theorem with $n=4$ and
$S=\{1,2,3\}$. It follows that \[ P(0,0,0) = \det \begin{bmatrix}
1&a_2&a_3&a_4\\
1& 1&a_2&a_3\\
0&1& 1&a_2\\
0&0&1&1
\end{bmatrix} = 1-3a_2+a_2^2+2a_3-a_4,\] as above. \label{newex43}
\end{example}
Our next result expresses the determinant in Theorem \ref{newthm41} as a
skew-Schur function of ribbon (also called rim-hook) type. Background on
skew-Schur functions can be found in \cite[Sec.~1.5]{mac}. For further
references and recent work, one can consult \cite{btw} or \cite{lp}. We
will use the elementary symmetric functions $e_r$ in countably many
variables defined as $e_0=1$ and \[ e_r = \sum_{1 \leq i_1 < i_2 < \dots <
i_r} x_{i_1} x_{i_2} \cdots x_{i_r}.\] From \cite[p.~20]{mac} these are
algebraically independent, so there is a homomorphism of the ring of
symmetric functions to $\mathbb{R}$ which sends each $e_i$ to $a_i$ from
Theorem \ref{newthm41}.

\begin{thm} \label{skewschur} With notation as in Theorem \ref{newthm41},
let $\lambda$ and $\mu$ be the partitions defined by
\[ \lambda_i=n-s_{i-1}-k+i-1 \ , \ \mu_i=n-s_{i}-k+i-1, \ \ 1 \leq i
\leq k+1.\] Let $\lambda',\mu'$ denote the transpose partitions of
$\lambda$ and $\mu$ and let $s_{\lambda'/\mu'}$ denote the corresponding
skew-Schur function, obtained by specializing the elementary symmetric
functions $e_i$ to equal $a_i$. Then
\[ P(t_1,\cdots,t_{n-1})= s_{\lambda'/\mu'}.\]
\end{thm}

\begin{proof} First note that the $\lambda_i,\mu_i$ defined above are all
non-negative. From Macdonald \cite[p.~71]{mac}, \begin{equation}
\label{jac} s_{\lambda'/\mu'} =
\det(e_{\lambda_i-\mu_j-i+j})_{i,j=1}^{k+1},\end{equation} and $e_r$ is
the $r$th elementary symmetric function. When each $e_r$ is specialized to
equal $a_r$, the quantity $s_{\lambda'/\mu'}$ becomes \[
\det(a_{\lambda_i-\mu_j-i+j})_{i,j=1}^{k+1} =
\det(a_{s_j-s_{i-1}})_{i,j=1}^{k+1},\] as desired. \end{proof}

One reason why Theorem \ref{skewschur} is interesting is that there are
non-obvious equalities between ribbon skew-Schur functions. The paper
\cite{btw} characterizes when two ribbon skew-Schur functions are equal;
analogous results for more general skew-Schur functions are in \cite{RSW}.
In particular, combining the results of \cite{btw} with Theorem
\ref{skewschur} one immediately obtains the fact that a stationary
one-dependent process is invariant under time reversal, i.e. \[
P(t_1,\cdots,t_n) = P(t_n,\cdots,t_1) \] for all $t_i$. For another proof
of invariance under time reversal, see \cite[p.~139]{pp}.

The determinantal formulae of Theorems \ref{newthm41} and \ref{skewschur}
suggest that a determinantal point process is lurking nearby. This is
indeed the case. In \ref{sect7} (see Corollary \ref{cor73}) we prove the
following.

\begin{cor}
A stationary one-dependent process as in Theorem \ref{newthm41} is
determinantal with kernel $K(x,y)=k(y-x)$ with
\begin{equation*}
\sum_{n\in\mathbb{Z}}k(n)z^n=-1/\sum_{j=1}^\infty a_jz^j.
\end{equation*}
\label{cor2}
\end{cor}

A remarkable development, connecting stationary one-dependent processes to
algebra appears in \citep[Chap.~7]{pp}. They consider a graded algebra
$A_0\oplus A_1\oplus A_2+\dots$ with $A_0=k$, a ground field, and
$A_iA_j\subseteq A_{i+j}$. For example, the space $k[x]$ of polynomials in
one variable has $A_0=k$, and $A_i$ spanned by $x^i$. They assume that
each $A_i$ is a finite dimensional vector space of dimension $\dim(A_i)$.
The algebra is \textit{quadratic} if $A_0=k$, the algebra is generated by
elements of $A_1$, and the relations defining the algebra are in $A_2$.
For example, the commutative polynomial ring $k[x_1,\dots,x_n]$ is
generated by $x_1,\dots,x_n$ and the quadratic relations
$x_ix_j-x_jx_i=0$. Note that algebras need not be commutative.

A technical growth condition on the $\dim(A_i)$ which we will not explain
here yields the Koszul algebras. These include \textit{many} natural
algebras occurring in mathematics (see \cite{fr} for a survey). As a
simple example, consider the commutative polynomial ring generated by
$x_1,x_2,\dots,x_n$ with the additional relations $x_ix_j=0$ for $(i,j)\in
E$, with $E$ the edge set of an undirected graph on $\{1,2,\dots,n\}$
(loops allowed). This is Koszul \citep[Chap.~2, Cor.~4.3]{pp}. Polishchuk
and Positselski \citep[Chap.~7, Cor.~4.3]{pp} prove the following
remarkable result.
\begin{thm} \label{deep}
To every Koszul algebra $A$ one can assign a stationary, one-dependent
process via
\begin{equation*}
P(X_1=X_2=\dots=X_{i-1}=1)=\frac{\dim(A_i)}{\dim(A_1)^i} \qquad\text{for
}i=1,2,\dots.
\end{equation*}
\label{newthm42}
\end{thm}
Part of the reason that Theorem \ref{deep} is substantive is that the
quantities $P(X_1=X_2=\dots=X_{i-1}=1)$ can not assume arbitrary values;
the determinants in Theorem \ref{newthm41} must be non-negative.

\begin{example}
Consider the commutative polynomial algebra generated by $x_1,x_2,\dots
x_b$ and the additional relations $x_i^2=0,\ 1\leq i\leq b$. The degree
$i$ part $A_i$ is spanned by square free monomials and has
$\dim(A_i)=\binom{b}{i},\ 0\leq i<\infty$. From Fact 2 of \ref{sec2}, we
see that the associated one-dependent process is precisely the carries
process of mod $b$ addition from our introduction. Note that since only
the $\dim(A_i)$ matter, the algebra in this example can be replaced by the
exterior algebra generated by $x_1,x_2,\dots, x_b$ and the relations
$x_ix_j=-x_jx_i$ for all $i,j$. \label{newex44}
\end{example}
\begin{example}
Consider the commutative polynomial algebra generated by $x_1,x_2,\dots,
x_b$. The degree $i$ part $A_i$ is spanned by monomials of degree $i$ and
has $\dim(A_i)=\binom{b+i-1}{i}$. From Example \ref{ex22} of \ref{sec2},
the associated one-dependent process is precisely the complement of the
carries process of mod $b$ addition. More generally, in \ref{sect7} we
show that the particle-hole involution of a (perhaps non-stationary)
one-dependent process with respect to any subset is one-dependent.
\label{newex45}
\end{example}
\begin{example} A PBW algebra is an algebra of the form \[
k[x_1,\cdots,x_n]/\langle x_ix_j: (i,j) \in S \subset [1,n]^2 \rangle,\]
where the variables $x_1,\cdots,x_n$ do not commute. Here it is useful to
think of $S$ as a directed graph with vertex set $\{1,\cdots,n\}$, loops
allowed. From \cite[p.~84]{pp}, PBW algebras are Koszul. As noted in
\citep[Chap.~7, Prop.~5.1]{pp}, the one-dependent process associated to a
PBW algebra can be described as follows: pick $U_1,U_2,\cdots$ i.i.d. in
$\{1,\cdots,n\}$ and let $X_i=h(U_i,U_{i+1})$, where $h(i,j)=0$ if and
only if $(i,j) \in S$. Indeed, the dimensions satisfy $\dim(A_0)=1,
\dim(A_1)=n$, $\dim(A_i)=$number of paths of length $i-1$ in the
complement of $S$, and so the chance of $i-1$ consecutive $1$'s in the
point process is equal to $\dim(A_i)/\dim(A_1)^i$.

The above argument shows that PBW algebras give rise to two-block-factor
processes. In fact they are dense in the set of all two-block factors
\cite[Chap.~7, Prop.~5.2]{pp}. For a generalization of these PBW
processes, as well as a simple argument that they are determinantal, see
Remark 5 after Theorem \ref{thm52}.
\end{example}

\begin{example}
Consider $2n$ points $x_1,x_2,\dots, x_{2n}$ in general position in
projective space $\mathbb{P}^n$. The coordinate ring of this projective
variety is known to be Koszul \citep[p.~42]{pp} with
$\dim(A_0)=1,\dim(A_1)=n+1,\dim(A_i)=2n$ for $i\geq2$. Thus
$P(X_1=X_2=\dots=X_{i-1}=1)=2n/(n+1)^i$ defines a one-dependent process.
When $n=1$, this is fair coin tossing. When $n=3$, it is the carries
process for multiplying quaternions introduced in \ref{sect6}.

For general $n$, these one-dependent processes can be described as
follows (and illustrates what one might consider to be a
``probabilistic'' description). Let $\{1,2,\cdots,n\}$ be ordered
cyclically and let $*$ be another symbol. Choose $U_i$ i.i.d. in
$\{*,1,\cdots,n\}$. Let $X_i=h(U_i,U_{i+1})$ where $h(*,x)=0$ for
all $x$, $h(x,*)=1=h(i,i+1)$ for $x \neq *$, and $h(i,j)=0$
otherwise. To see that this works, note that $P(X_1=1)=2n/(n+1)^2$,
since one has to choose $x$ different from $*$ first and then $*$ or
$x+1$ next. Similarly, for the chance of $i-1$ $1$'s in a row, the
first choice can be anything but $*$, then the next $i-2$ choices
are determined and the last can be one of two.

\label{newex46}
\end{example}

We are certain that natural Koszul algebras will lead to natural point
processes. We note further that \citep[Sec.~7.6]{pp} shows how natural
operations on algebras preserve the Koszul property. These include the
operations of union and complements that we work with in \ref{sect7}.

\section{Descents in permutations}\label{sect5}

A permutation $\sigma\in S_n$, the symmetric group, has a descent at $i$
if $\sigma(i)>\sigma(i+1)$. The set of such $i$ forms the descent set
$D(\sigma)$. If the base $b$ in previous sections is large and $n$ stays
fixed, a string of $n$ digits will have no repeated values and the descent
theory of \ref{sec2} and \ref{sec3} above becomes descent theory for
random permutations. This is a venerable subject. Stanley
\cite{sta1,sta2.1} reviews the basics. Some modern highlights are
Solomon's descent algebra \cite{sol}, the connections with the free Lie
algebra \cite{gars,reu}, Gessel's theory of enumerating permutations by
descents and cycle structure \cite{pd77,ges}, quasi-symmetric functions
and the theory of riffle shuffling \cite{ful02,sta3}. Any Coxeter group
has its own descent theory \cite{sol}. There is some indication that
arithmetic carries can be carried over \cite{pd173}.

This section introduces three examples where the descent set can be shown
to be a one-dependent determinantal point process: uniform choice of
permutations, non-uniform choice from the Mallows model, and independent
trials. Then it shows that the closely related notion of the connectivity
set of a permutation also yields a determinantal process. Finally,
connections with binomial posets are mentioned.
\begin{example} \label{des1} \textbf{The descent set of a uniformly chosen
    permutation.}\quad Consider the formula of \ref{sec2} Fact 5 for
  the chance that a random $b$-ary string of length $n$ has descents
  exactly at $S\subseteq[n-1]$. Passing to the limit, as
  $b\nearrow\infty$ using $\binom{a+b-1}{b-1}\sim\frac{b^a}{a!}$ gives
  a classical formula for the chance that a uniformly chosen
  $\sigma\in S_n$ has descents exactly at $S$:
\begin{equation}
P_n\left(D(\sigma)=S\right)=\text{det}\left[1\big/\left(s_{j+1}-s_i\right)!\right].
\label{41}
\end{equation}
For $S=1\leq s_1<s_2<\dots<s_k\leq n-1$, the determinant is of a
$(k+1)\times(k+1)$ matrix with $(i,j)$ entry $1/(s_{j+1}-s_i)!,
(i,j)\in[0,k]\times[0,k]$ and $s_0=0, s_{k+1}=n$. This formula is
originally due to MacMahon \cite{mcm}. See Stanley \citep[p.~69]{sta2} or
Gessel and Viennot \cite{gv} for modern proofs. \label{exam41}
\end{example}

Equation \eqref{41} allows us to see that the descent set of a uniformly
random permutation, treated as a point process, is determinantal. The
proof follows from \eqref{41} and Corollary \ref{cor1} in \ref{sect7}.
\begin{thm}
There exists a stationary point process on $\mathbb{Z}$, call it $P$, such
that its restriction to any interval of length $n-1$ coincides with $P_n$
of \eqref{41}. The process $P$ is one-dependent and determinantal. Its
correlation kernel $K(x,y)=k(y-x)$ with
\begin{equation}
\sum_{m\in\mathbb{Z}}k(m)z^m=\dfrac1{1-e^{z}}. \label{42}
\end{equation}
\label{thm41}
\end{thm}
\begin{rems}\

\begin{enumerate}
\item Theorem \ref{thm41} can also be seen by passing to the limit in
  Theorem \ref{thm31}. Replacing $t$ in the generating function
  $1/(1-(1-t)^b)$ by $-z/b$ and letting $b\nearrow\infty$ gives
  \eqref{42}.

\item The Bernoulli numbers $B_n$ are defined by
  $z/(e^z-1)=\sum_{n=0}^\infty B_nz^n/n!$. These are very
  well-studied. It is known that $B_{2i+1}=0$ for $i\geq1$, $B_0=1,
  B_1=-\tfrac12, B_2=\tfrac16, B_4=-\tfrac1{30}, B_6=\tfrac1{42},
  \dots$. We see that $k(m)=-B_{m+1}/(m+1)!$.

\item From these calculations, the kernel $K$ is
\begin{equation*}
K=\begin{bmatrix}
\tfrac12&-\tfrac1{12}&0           &\tfrac1{720}&\cdots\\
-1      &\tfrac12    &-\tfrac1{12}&0           &\ddots\\
        &-1          &\ddots      &\ddots      &\ddots\\
        &            &\ddots      &\ddots      &\ddots\\
        &            &            &\ddots      &\ddots
\end{bmatrix}\ .
\end{equation*}

\item The correlation functions of $P$ (and $P_n$) can be computed
  explicitly by passage to the limit from Fact 5 of \ref{sec2}. Let
  $A=A_1\cup A_2\cup\dots\cup A_k$ be a decomposition of the finite
  set $A\subseteq \mathbb{Z}$ into disjoint blocks of adjacent
  integers. If $|A_i|=a_i$,
\begin{equation*}
\rho(A)=\prod_{i=1}^k\dfrac1{(a_i+1)!}.
\end{equation*}

\item For any $n=1,2,\dots$, let $K_n(i,j)=k(j-i)$ be the $n\times n$
  top left corner block of $K$. Let $d(\pi)$ denote the number of descents
  of $\pi$ and let \[ A_{n+1} = \sum_{\pi \in S_{n+1}} x^{d(\pi)} \] be
  the $(n+1)$st Eulerian polynomial (see e.g., \citep[p.~22]{sta2}). Then
 item (c) preceding Theorem \ref{thm31} yields that \begin{equation*}
A_{n+1}(x)=\text{det}\left(I+(x-1)K_n\right)(n+1)!
\end{equation*} It is known \citep[p.~292]{com} that $A_{n+1}(x)$ has all real
zeros $\alpha_1,\alpha_2,\dots,\alpha_n$. From the development in
\ref{sec3}, $N_{n+1}$, the number of descents in a random permutation from
$S_{n+1}$, is the sum of $n$ independent Bernoulli random variables with
success rates $(1-\alpha_j)^{-1},\ 1\leq j\leq n$. One easily computes
\begin{align*}
E(N_{n+1})=\text{tr}(K_n)=\dfrac{n}2,\qquad\text{Var}(N_{n+1})&=\text{tr}
\left(K_n-K_n^2\right)\\
&=\dfrac{n}2-\left(\frac{5n}{12}-\dfrac16\right)= \frac{n+2}{12}.
\end{align*}
It follows by Harper's method (\cite{har},\cite{Pi}) that one has a
central limit theorem:
\begin{equation*}
\frac{N_{n+1}-\tfrac{n}2}{\sqrt{\frac{n+2}{12}}}\Longrightarrow\mathcal{N}(0,1).
\end{equation*}

\item Chebikin \cite{ch} studied the ``alternating descent set'' $A(\sigma)$ of a
random permutation, defined as the set of positions at which the
permutation has an {\it alternating descent}, that is an ordinary descent
if $i$ is odd, or an ascent if $i$ is even. Combining Lemma 2.3.1 of
\cite{ch} with the argument used to prove \eqref{41}, one obtains that for
$|S|=k$,
\[ P_n(A(\sigma)=S) = \det[E_{s_{j+1}-s_i}/(s_{j+1}-s_i)!]_{i,j=0}^k,\] where $E_n$
is the nth Euler number defined by $\sum_{n \geq 0} E_n z^n/n! = \tan(z) +
\sec(z),$ and $s_0=0, s_{k+1}=n$. Applying our Theorem \ref{thm51} in
\ref{sect7}, it follows that that $P_n$ is obtained by restricting to any
interval of length $n-1$ the stationary, one-dependent, determinantal
process with correlation kernel $K(x,y)=k(y-x)$ with \[ \sum_{m \in
\mathbb{Z}} k(m) z^m= \frac{1}{1-(\tan(z)+\sec(z))^{-1}}.\]

\item Modern combinatorics suggests a host of potential
  generalizations of Theorem \ref{thm41}. Let $P$ be a partial order on $[n]$.
  A linear extension of $P$ is a permutation $\sigma\in S_n$ such that if $i$
  is less than $j$ according to $P$, then $\sigma(i)<\sigma(j)$. Let
  $\mathcal{L}(P)$ denote the set of linear extensions of $P$. For each
  $\sigma\in\mathcal{L}(P)$, let
\begin{align*}
D(\sigma)&=\{i:\sigma(i)>\sigma(i+1)\},&\text{the
  descent set of }\sigma,\\
\wedge(\sigma)&=\{i:\sigma(i-1)<\sigma(i)>\sigma(i+1)\},&\text{the
  peak set of }\sigma.
\end{align*}
Choosing $\sigma$ uniformly in $\mathcal{L}(P)$ gives two point processes.
This section has focused on ordinary descents, that is, $P$ is the trivial
poset with no restrictions. There are many indications that ``natural''
posets will give rise to determinantal point processes. For background and
first results, see Brenti \cite{bren} or Stembridge \cite{stem}.

\item Finally, we mention that it would be interesting to find analogs of
Theorem \ref{thm41} for other Coxeter groups. For the hyperoctahedral
group $B_n$ consisting of the $2^nn!$ signed permutations, descents can be
defined using the linear ordering \[ 1<2<\dots<n<-n<\dots<-2<-1.\] Say
that
\begin{enumerate}
\item $\sigma$ has a descent at position $i$ $(1 \leq i \leq n-1)$ if
$\sigma(i) > \sigma(i+1)$.
\item $\sigma$ has a descent at position $n$ if $\sigma(n)<0$.
\end{enumerate}
For example, $-4 \ -1 \ 3 \ 2 \ 5 \ \in B_3$ has descent set $\{2,3,5\}$.
Reiner \cite{Re} shows that the chance that a random element of $B_n$ has
descent set $\{s_1,\cdots,s_r\} \subseteq \{1,2,\cdots,n\}$ is
$\det(a(i,j))_{i,j=0}^r$ where
\[ a(i,j) = \left\{ \begin{array}{ll}
\frac{1}{(s_{j+1}-s_i)!} & \mbox{if $0 \leq j \leq r-1$}\\
\frac{1}{2^{n-s_i}(n-s_i)!} & \mbox{if $j=r$}
\end{array} \right.\] Here $s_0=0$. Theorem \ref{thm51} in \ref{sect7}
shows that the resulting processes is determinantal (but not stationary).
\end{enumerate}
\end{rems}
\begin{example} \label{des2} \textbf{Descents from non-uniform distributions on
    permutations.}\quad In a variety of applications,
  \textit{non}-uniform distributions are used on permutations. A
  widely used model is
\begin{equation}
P_\theta(\sigma)=Z^{-1}\theta^{d(\sigma,\sigma_0)},\qquad\sigma\in S_n.
\label{43}
\end{equation}
Here, $0<\theta\leq1$ is a fixed concentration parameter, $\sigma_0\in
S_n$ is a fixed location parameter, $d(\sigma,\sigma_0)$ is a metric on
$S_n$ and $Z=\sum_\sigma\theta^{d(\sigma,\sigma_0)}$ is the normalizing
constant. These are called ``Mallows models through the metric $d$.'' When
$\theta=1$, $P_1$ is the uniform distribution. For $0<\theta<1$,
$P_\theta$ is largest at $\sigma_0$ and falls off as $\sigma$ moves away
from $\sigma_0$. See \cite{cri}, \cite{dr} or \cite{mar} for background
and references. \label{exam42}
\end{example}

Perhaps the most widely used metric is
\begin{equation}\begin{aligned}
d(\sigma,\sigma_0)&=\text{minimum number of pairwise adjacent transpositions}\\
&\qquad\text{to bring $\sigma$ to $\sigma_0$},\\
&=I\left(\sigma\sigma_0^{-1}\right),\text{ the number of inversions in }\sigma
\sigma_0^{-1}.
\end{aligned}\label{44}
\end{equation}
This is called ``Kendall's tau'' in the statistics literature. References,
extensions and properties are in \cite[Sect.~4]{d}, \cite{flig},
\cite{mal}. We mention that the normalizing constant of the Mallows model
through Kendall's tau is given by $Z=\prod_{i=1}^n
\frac{\theta^i-1}{\theta-1}$ \cite[p.~21]{sta2}.

With $\theta,\sigma_0$ and $d$ fixed, one may ask any of the usual
questions of applied probability and enumerative combinatorics: ``Picking
a permutation randomly from $P_\theta(\cdot)$, what is the distribution of
the cycle structure, longest increasing subsequence, \textellipsis?'' For
general metrics, these are open research problems. Even the algorithmic
task of sampling from $P_\theta$ leads to difficult problems. See
\cite{dh} or \cite{d}.

For Mallows model $P_\theta$ through Kendall's tau \eqref{44} and
$\sigma_0=\text{id}$, we show that the descent set of $\sigma,\
D(\sigma)=\{i:\sigma(i+1)<\sigma(i)\}$ forms a determinantal point process
with simple properties. In what follows we use the $q$-analog notation
that $i_q=\frac{q^i-1}{q-1}$ and $n!_q=\prod_{i=1}^n i_q$, and take
$\theta=q$ as is conventional in combinatorial work in the subject.
\begin{prop}
Let $P_q$ $(0<q<1)$ be Mallows model through Kendall's tau \eqref{43},
\eqref{44} on $S_n$.
\begin{itemize}
\item[(a)] The chance that a random permutation chosen from $P_\theta$ has
  descent set $s_1<s_2<\dots<s_k$ is
\begin{equation*}
\emph{det}\left[\frac{1}{(s_{j+1}-s_i)!_q}\right]_{i,j=0}^k,
\end{equation*}
with $s_0=0$ and $s_{k+1}=n$.

\item[(b)] The point process associated to $D(\sigma)$ is stationary,
one-dependent, and determinantal with kernel $K(x,y)=k(y-x)$, where
\begin{equation*}
\sum_{m \in \mathbb{Z}} k(m)z^m=\dfrac1{1-\left(\sum_{m=0}^\infty
z^m\big/m!_q \right)^{-1}}.
\end{equation*}
\item[(c)] The chance of finding $k$ descents in a row is
$q^{\binom{k+1}{2}}/(k+1)!_q$.
  In particular, the chance of a descent at any given location is
  $q/(q+1)$. The number of descents has mean
  $\mu(n,q)=\frac{(n-1)q}{q+1}$ and variance
  $q \left[ \frac{(q^2-q+1)n-q^2+3q-1}{(q^2+q+1)(q+1)^2} \right]$.
  Normalized by its mean and variance, the number of descents has a
  limiting standard normal distribution.
\end{itemize}
\label{prop1}
\end{prop}
\begin{proof}
  Part (a) of Proposition \ref{prop1} follows from Stanley
  \cite[Ex.~2.2.5]{sta2}. Part (b) follows from Corollary \ref{cor1} in
  \ref{sect7} and elementary calculations. To compute the chance of $k$
  consecutive descents, note by stationarity that this is the chance of
  the permutation $\sigma(i)=k+2-i$ in $S_{k+1}$ under the $P_q$ measure
  on $S_{k+1}$. This probability is $q^{\binom{k+1}{2}}/(k+1)!_q$, as
  needed. One calculates from part (b) that \[ k(-1)=1, \
  k(0)=\frac{q}{q+1}, \ k(1)=\frac{q^2}{(q+1)^2(q^2+q+1)}.\] Recall
  from the remarks preceding Theorem \ref{thm31} that \[ E(x^N)=\text{det}\left(I+(x-1)
  K\right),\] where $N$ is the total number of particles. Thus
\[ E(|D(\sigma)|) = tr(K_{n-1}) = (n-1)k(0), \]
\[ Var(|D(\sigma)|) = tr(K_{n-1}-K_{n-1}^2) = (n-1)k(0) - (n-1)k(0)^2
- 2(n-2)k(1)k(-1),\] and the proof of part (c) is complete.
\end{proof}

{\it Remarks.} \begin{enumerate} \item
  From part (c), the correlation functions are given explicitly, just
  as in remark 4 following Theorem \ref{thm41}.

  \item It would be interesting to have a direct combinatorial proof of
  the one-dependence in Proposition \ref{prop1}.
\end{enumerate}
\begin{example} \textbf{Descents for independent trials}\quad
The carries in the carries process have an equivalent description in
terms of descents in an independent and uniform sequence. In this
example we generalize to descents in a sequence $Y_1,Y_2,\dots,Y_n$
with $Y_i$ independent and identically distributed with
$P(Y_1=j)=p_j$. Here $0\leq j\leq b-1$ to match previous sections
(and of course $0\leq p_j\leq1, \sum p_j=1$), but in fact $b=\infty$
can be allowed. The descents are related to the order statistics for
discrete random variables \cite{abn,eld} but we have not seen this
problem previously treated in the probability or statistics
literature. Of course, if $Y_i$ are independent and identically
distributed from a continuous distribution, the descent theory is
the same as the descent theory for a random permutation.
\label{exam54}
\end{example}

It is useful to have the complete homogeneous symmetric functions
\begin{equation*}
h_k(x_1,x_2,\dots,x_m)=\sum_{1\leq i_1\leq\dots\leq i_k\leq
  m}x_{i_1}x_{i_2}\dots x_{i_k}.
\end{equation*}
\begin{thm}
If $Y_1,Y_2,\dots,Y_n$ are independent with common probability
$P(Y_1=j)=p_j$, $0\leq j\leq b-1$, the descents in $Y_1,Y_2,\dots,Y_n$
form a stationary one-dependent point process on $[n-1]$ with
\begin{equation*}
P_n(S)=\emph{det}\left(h_{s_{j+1}-s_i}\right)_{0\leq i,j\leq
    k}\qquad\text{with }s_0=0,s_{k+1}=n,
\end{equation*}
for $S=\{s_1<s_2<\dots<s_k\}$. The $h_j$ are evaluated at
$p_0,\dots,p_{b-1}$. Moreover, this process is determinantal with
$K(i,j)=k(j-i)$ and $\sum_{l=-\infty}^\infty
k(l)z^l=\tfrac1{1-\prod_{i=0}^{b-1}(1-p_iz)}$. \label{thm52}
\end{thm}
\begin{proof}
Let $\alpha_n(S)=P_n(\text{des}(Y_1,Y_2,\dots,Y_n)\subseteq
S)=h_{s_1}h_{s_2-s_1}\dots h_{n-s_k}$. The determinant formula follows
from the argument for Fact 5 in \ref{sec2}. Now, Corollary \ref{cor1}
implies that $\hat{k}(z)=\tfrac1{1-\hat{e}(z)^{-1}}$, with
$\hat{e}(z)=\sum_{l=0}^\infty h_l z^l=\prod_{i=0}^{b-1}(1-p_iz)^{-1}$ from
e.g., Macdonald \cite[p.~21]{mac}.
\end{proof}
\begin{rems}\

\begin{enumerate}
\item Setting $p_i=1/b,\ \hat{k}(z)=\tfrac1{1-(1-z/b)^b}$. This is
  equivalent to our result in Theorem \ref{thm31}, by changing
  variable $z$ to $z/b\to t$, which doesn't change correlations.

\item From the theorem, the associated point process is
  one-dependent. If $X_i=1$ or $0$ as there is a descent at $i$ or
  not, then clearly
\[ P(X_1=X_2=\dots=X_i=1) = e_{i+1}(p_0,\cdots,p_{b-1}),\] where
$e_r(x_1,\cdots,x_m) = \sum_{1 \leq i_1<\cdots<i_r \leq m} x_{i_1} x_{i_2}
\cdots x_{i_m}$ is the $r$th elementary symmetric function. For example,
\begin{align*}
P(X_i=1)&=\dfrac12-\dfrac12\sum_ip_i^2,\\
P(X_i=X_{i+1}=1)&=\dfrac16-\dfrac12\sum_ip_i^2 + \dfrac13\sum_ip_i^3.
\end{align*}
This allows a straightforward formulation (and proof) of the central
limit theorem for the number of descents. As in \ref{sec2}, simple
expressions for the correlation functions are available in terms of
$P(X_1=X_2=\dots=X_i=1)$.

\item As we show in Corollary \ref{skew2}, one can express $P_n(S)$
as a skew-Schur function of ribbon type. For special choices of the
$p_i$'s, this symmetric function interpretation may facilitate
computation. For example, \cite[Sect.~1.5]{mac} gives
$s_{\lambda/\mu}=\sum_{\nu} c_{\mu \nu}^\lambda s_{\nu}$ with $c_{\mu
\nu}^\lambda$ the Littlewood--Richardson coefficients and $s_{\nu}$ the
Schur function. If $p_i=zq^i$ ($i=0,1,\cdots,b-1$, $z$ a normalizing
constant), then $s_{\nu}$ is non-0 only if $\nu$ has length $\leq b$, in
which case
\begin{equation*}
s_{\nu}=z^{|\nu|}q^{n(\nu)}\prod_{x\in
\nu}\frac{1-q^{b+c(x)}}{1-q^{h(x)}};
\end{equation*}
see Macdonald \cite[p.~44]{mac} for an explanation of notation.

\item
Coming back to carries, if a column of numbers is chosen
  independently on $\{0,1,\dots,b-1\}$ from common distribution
  $\{p_i\}_{i=0}^{b-1}$, the remainder column after addition evolves
  as a Markov chain with starting distribution $\{p_i\}$ and
  transition matrix $P(i,j)=p_{j-i}$ (indices mod $b$). This is a
  circulant with known eigenvalues and eigenvectors. We do not know
  the distribution of carries.

\item There is a sweeping generalization of this example which leads
 to a large collection of determinantal point processes. Let $R$ be an
 arbitrary subset of $[N]\times[N]$. Fix a probability distribution
  $\theta_1,\theta_2,\dots,\theta_N$ on $[N]$. A point process $P_n$ on
  $[n-1]$ arises as follows: pick $Y_1,Y_2,\dots,Y_n$ independently
  from $\{\theta_i\}$. Let $X_i={\footnotesize\begin{cases}1&\text{if
      }(Y_i,Y_{i+1}) \notin R\\0&\text{otherwise}\end{cases}},\ 1\leq
  i\leq n-1$. Let $S=\{i:X_i=1\}$. Using essentially the arguments
  above, it is proved in Brenti \cite[Th.~7.32]{bren} that
\begin{equation*}
P_n(S)=\text{det}\left[h_{s_{j+1}-s_i}^R\right]_{i,j=0}^k\qquad\text{for
}S=\{s_1<\dots<s_k\}, s_0=0, s_{k+1}=n,
\end{equation*}
where $h_j^R=h_j^R(\theta_1,\dots,\theta_N):=P_j(\emptyset)$. Here as
usual, $h_0^R=1,h_j^R=0$ for $j<0$. This falls into the domain of Theorem
\ref{thm51} and its corollary. It follows that $P_n$ is determinantal with
kernel $K(i,j)=k(j-i),\ \sum_{l \in \mathbb{Z}}
k(l)z^l=\tfrac1{\left(1-\left(\sum_{j=0}^\infty
h_j^Rz^j\right)^{-1}\right)}$.

Our examples are the special case
$R(i,j)={\footnotesize\begin{cases}1&\text{if }i\leq
    j\\0&\text{otherwise}\end{cases}}$. We are certain that there are many
other interesting cases.
\end{enumerate}
\end{rems}

\begin{example} \label{connect} \textbf{Connectivity set of a permutation.}\quad

Stanley \cite{sta4} studied a ``dual'' notion to the descent set of a
permutation, which he called the connectivity set. For $\sigma \in S_n$,
the connectivity set $C(\sigma)$ is defined as the set of $i$, $1 \leq i
\leq n-1$, such that $\sigma \in S_{1,\cdots,i} \times S_{i+1,\cdots,n}$.
For example, the permutation $\sigma= 3 \ 1 \ 2 \ 5 \ 4 \ 6$ satisfies
$C(\sigma)=\{3,5\}$. The connectivity set also arises in the analysis of
quicksort, where it is called the set of splitters of a permutation
\cite[Sec.~2.2]{wil}.

We will prove that the connectivity set of a random element of $S_n$ is a
determinantal point process and determine its correlation kernel. For this
the following two facts are helpful:

\begin{enumerate}
\item  A permutation $\sigma \in S_n$ is called {\it indecomposable} or {\it
connected} if $C(\sigma)=\emptyset$. Comtet \cite[Exer.~VII.16]{com} shows
that the number $f(n)$ of connected permutations in $S_n$ satisfies
\begin{equation} \sum_{n \geq 1} f(n)x^n = 1 - \frac{1}{\sum_{n \geq 0} n!x^n}.
\label{genf} \end{equation} (An asymptotic expansion of $f(n)$ is given in
\cite[Exer.~VII.17]{com}).

\item \cite[Prop.~1.1]{sta4} Letting $S=\{s_1<s_2<\cdots<s_k\}$ be a
subset of $\{1,2,\cdots,n-1\}$, \begin{equation} \label{cont} \left|
\{\sigma \in S_n: S \subseteq C(\sigma) \} \right| = s_1! (s_2-s_1)!
\cdots (s_k-s_{k-1})! (n-s_k)!.\end{equation}
\end{enumerate}

The following is our main result.

\begin{thm} Let $C(\sigma)$ denote the connectivity set of a permutation
$\sigma$ chosen uniformly at random from $S_n$. The point process
corresponding to $C(\sigma) \cup \{0,n\}$ is determinantal, with state
space $\{0,\cdots,n\}$ and correlation kernel $K(x,y)$ satisfying
\[
\begin{array}{ll}
K(0,y)= 1 & \mbox{all $y$}\\
K(n,y)= \delta_{n,y} & \mbox{all $y$} \\
K(x,y)= \frac{1}{{n \choose x}} & \mbox{$0,n \neq x \leq y$}\\
K(x,y)= \frac{1}{{n \choose x}} - \frac{1}{{n-y \choose n-x}} & \mbox{$0,n
\neq x>y$}
\end{array} \]
\end{thm}

Note that in the statement of the theorem, $0,n$ are always points of the
process.

\begin{proof} The first step is to observe that $C(\sigma) \cup \{0,n\}$
can be obtained as a trajectory of a certain Markov chain, started at $0$,
with transition probabilities \begin{equation} P(i,j) =
\frac{(n-j)!f(j-i)}{(n-i)!}. \label{trans} \end{equation} Here we take
$f(0)=0$, so $P(i,i)=0$ for all $i$. To prove \eqref{trans}, note that
$C(\sigma)=\{i_1,\cdots,i_k\}$ if and only if the following events
$E_1,\cdots,E_{k+1}$ occur:
\begin{itemize}
\item $E_1$: $\{\pi(1),\cdots,\pi(i_1)\}=\{1,\cdots,i_1\}$ and $\pi$ restricted to
$\{1,\cdots,i_1\}$ is indecomposable.
\item $E_2$: $\{\pi(i_1+1),\cdots,\pi(i_2)\}=\{i_1+1,\cdots,i_2\}$ and $\pi$
restricted to $\{i_1+1,\cdots,i_2\}$ is indecomposable.
\item $\cdots$
\item $E_{k+1}$: $\{\pi(i_k+1),\cdots,\pi(n)\}=\{i_k+1,\cdots,n\}$ and $\pi$
restricted to $\{i_k+1,\cdots,n\}$ is indecomposable.
\end{itemize} Letting $f(n)$ denote the number of indecomposable permutations
in $S_n$, the probability of $E_1$ is clearly $\frac{f(i_1)(n-i_1)!}{n!}$.
The probability of $E_2$ given $E_1$ is
$\frac{f(i_2-i_1)(n-i_2)!}{(n-i_1)!}$ and the probability of $E_3$ given
$E_1,E_2$ is $\frac{f(i_3-i_2)(n-i_3)!}{(n-i_2)!}$, etc., as claimed.

Now \cite[Thm.~1.1]{bo} implies that the point process on $C(\sigma) \cup
\{0,n\}$ is determinantal with correlation kernel \begin{equation}
\label{ck} K(x,y) = \delta_{0,x} + Q(0,x) - Q(y,x) \end{equation} where \[
Q=P+P^2+P^3+ \cdots.\] To compute $Q$, let $[x^r] g(x)$ denote the
coefficient of $x^r$ in a series $g(x)$, and
note that \begin{eqnarray*} P^2(i,j) & = & \sum_l P(i,l)P(l,j)\\
& = & \sum_l \frac{(n-l)!f(l-i)}{(n-i)!} \frac{(n-j)!f(j-l)}{(n-l)!}\\
& = & \frac{(n-j)!}{(n-i)!} [x^{j-i}] \left(1 - \frac{1}{\sum_{n \geq 0}
n!x^n} \right)^2. \end{eqnarray*} The last equality used \eqref{genf}. A
similar computation shows that
\[ P^r(i,j) = \frac{(n-j)!}{(n-i)!} [x^{j-i}] \left(1 - \frac{1}{\sum_{n
\geq 0} n!x^n} \right)^r,\] and thus
\begin{eqnarray*}
Q(i,j) & = & \frac{(n-j)!}{(n-i)!} [x^{j-i}] \frac{\left(1 -
\frac{1}{\sum_{n \geq 0} n!x^n} \right)}{\left(\sum_{n \geq 0}
n!x^n \right)^{-1}}\\
& = & \frac{(n-j)!}{(n-i)!} [x^{j-i}]\left(\sum_{n \geq 1} n!x^n
\right).
\end{eqnarray*} Thus $Q(i,j)=1/{n-i \choose n-j}$ if $i<j$, and is
$0$ otherwise. The theorem follows from this and \eqref{ck}.
\end{proof}

{\it Remarks.}
\begin{enumerate}
\item The connectivity set $C(\sigma)$ is a simple example of a
determinantal process which is not one-dependent. Indeed, from
\eqref{cont}, $P(1 \in C(\sigma))=\frac{1!(n-1)!}{n!}$, $P(3 \in
C(\sigma))=\frac{3!(n-3)!}{n!}$, and $P(1,3 \in C(\sigma)) =
\frac{1!2!(n-3)!}{n!} \neq P(1 \in C(\sigma)) P(3 \in C(\sigma))$.

\item Unlike the other point processes considered in this paper, the
expected number of points tends to $0$ as $n \rightarrow \infty$.
Indeed, applying \eqref{cont} gives that \[ E|C(\sigma)| =
\sum_{i=1}^{n-1} P(i \in C(\sigma)) = \sum_{i=1}^{n-1} \frac{1}{{n
\choose i}} \rightarrow 0 \] as $n \rightarrow \infty$.
\end{enumerate}

\end{example}

\subsection*{Binomial posets}

In Doubilet--Rota--Stanley \cite{drs}, binomial posets were introduced as
a unifying mechanism to ``explain'' the many forms of generating functions
that appear in enumerative combinatorics. Briefly, $\mathcal{X}$ is a
binomial poset if for every interval $[x,y]$, the length of all maximal
chains is the same (say $n_{xy}$) and the number of maximal chains in an
interval of length $n$ does not depend on the particular interval. This
number $B(n)$ is called the factorial function. For the usual Boolean
algebra of subsets, $B(n)=n!$. For the subspaces of a vector space over a
finite field, $B(n)=n!_q$ ($q$-factorial). There are many further
examples.

Stanley \cite{sta1} has shown that many enumerative formulae generalize
neatly to the setting of binomial posets. Some of these developments give
new determinantal point processes. Here is an example.
\begin{example} \label{int} (union of descent sets)

\begin{thm}
Pick $\sigma_1,\sigma_2,\dots,\sigma_n$ independently from Mallows model
through Kendall's tau \eqref{43} on the symmetric group $S_n$. Let
$S=\{s_1<\cdots<s_k\}$ be the union of their descent sets. Then
\begin{equation*}
P_n(S)=\emph{det}\left[\dfrac1{\left(\left(s_{j+1}-s_i\right)!_q\right)^r}\right],
\end{equation*} where $s_0=0, s_{k+1}=n$. The associated point process is
stationary, one-dependent, and determinantal with $K(x,y)=k(y-x)$ for
$\hat{k}(z)=\tfrac1{1-1/\hat{e}(z)}$ where $\hat{e}(z)=\sum_{l=0}^\infty
z^l/(l!_q)^r$. \label{thm53}
\end{thm}

\begin{proof} The formula for $P_n(S)$ is Corollary 3.2 of \cite{sta1},
so the theorem follows from one of our main results proved in \ref{sect7}
(Corollary \ref{cor1}).
\end{proof}

When $q=1$ and $r=2,\ \hat{e}(z)=I_0(2 \sqrt{z})$, the classical
modified Bessel function. Feller \cite[\S II.7]{fel} develops a host
of connections with stochastic processes.
\begin {rem}
  The intersection of two independent point processes with correlation
  functions $\rho^1(A),\rho^2(A)$ is a point process with correlation
  function $\rho^1(A)\rho^2(A)$. If both processes are determinantal,
  this is a product of determinants, but in general the intersection or union
  of determinantal processes is not determinantal. In Subsection \ref{all},
  we show that all one-dependent processes are determinantal. Since the
  intersection of descent sets of independent permutations is
  one-dependent, and taking complements preserves one-dependence,
  this gives another proof that the union process in Theorem \ref{thm53} is
  determinantal.
\end{rem}
\end{example}

We have not pursued other examples but again believe there is much
else to be discovered.

\section{More general carries}\label{sect6}

Consider the quaternions $Q_8=\{\pm1,\pm i,\pm j,\pm k\}$. The
center $Z$ of $Q_8$ is $\{\pm1\}$. Choose coset representatives
$X=\{1,i,j,k\}$ for $Z$ in $Q_8$ so any element can be uniquely
represented as $g=zx$. The coset representatives are multiplied by
the familiar rule $k\curvearrowright i\curvearrowright j
\curvearrowright k$ so $ij=k$ and $kj=-i$, etc. If we multiply
$g_1g_2\dots g_k=(z_1x_1)(z_2x_2)\dots(z_kx_k)$, the
$z_i\in\{\pm1\}$ can all be moved to the left and we must multiply
$x_1 x_2\cdot\dots\cdot x_k$, keeping track of the ``carries'', here
$\pm1$. Evidently, if $\{g_i\}$ are chosen uniformly at random in
$Q_8$, both $\{z_i\}$ and $\{x_i\}$ are independent and uniform in
$Z$ and $X$. Thus we have the following problem: choose
$X_1,X_2,\dots,X_k$ uniformly in $X$ and multiply as $X_1$,
$X_1X_2$, $(X_1X_2)X_3$, \dots,
\begin{center}
$X_1$\\
$X_2$\\
$X_3$\\
$\vdots$\\
$X_k$
\end{center}
This gives a process of remainders and carries as in \ref{sec1}.
\begin{example}
\begin{equation*}\begin{array}{rrr}
k&\cdot&k\\
k&&1\\
i&\cdot&i\\
i&&1\\
k&\cdot&k\\
j&\cdot&i\\
k&\cdot&j\\
i&&k
\end{array}\qquad\text{gives }(-1)^5k=-k.
\end{equation*}
It is almost obvious that the carries form a stationary, one-dependent,
two-block process with $P(i-1$ ones in a row) $=\tfrac{6}{4^{i}},\ 2\leq
i<\infty$. Further, the remainders in the second column are independent
and uniform on $\{1,i,j,k\}$ and there is a simple ``descent'' rule which
determines the joint law of the dots (Example \ref{ex72} below).

We also mention, by comparison with Example \ref{newex46} in
\ref{newsec4}, that the carries process for the quaternions is the same
point process that arises from the coordinate ring of 6 generic points in
projective space $\mathbb{P}^3$. \label{ex71}
\end{example}

One natural generalization where all goes through is to consider a finite
group $G$ and  a normal subgroup $N$ contained in the center of the group.
The factor group $F=G/N$ has elements labeled $1,\sigma,\tau,\dots$. We
may choose coset representatives $t(1)=1,t(\sigma),t(\tau),\dots$ and any
$g=nt(\sigma)$, uniquely. While sometimes $t(\sigma)t(\tau)=t(\sigma
\tau)$, in many cases this fails; but we may choose correction factors
$f(\sigma,\tau)$ in $N$ (often called a ``factor set'') so that
$t(\sigma)t(\tau)=t(\sigma \tau)f(\sigma,\tau)$. Once $t(\sigma)$ are
chosen, the $f(\sigma,\tau)$ are forced.
\begin{example}
If $G=C_{100},\ N=C_{10}$ (thought of as a subgroup
$\{0,10,20,\dots,90\}$), the natural choice of coset representatives for
$G/N \cong \{0,1,\dots,9\}$ is $t(i)=i\in G$. Of course,
\begin{equation*}
t(i)+t(j)=t(i+j)+f(i,j),\qquad\text{with }f(i,j)=\begin{cases}
1&\text{if }i+j\geq10\\
0&\text{if }0\leq i+j<10.\end{cases}
\end{equation*} It is natural to ask if the choice of cosets matters. To see that it can,
consider $C_{10}$ in $C_{100}$ and choose coset representatives as
$0,11,22,43,44,45,46,47,48,49$. The sum of $11$ and $22$ requires a carry
of $90$.

A lovely exposition of carries as cocycles is in \cite{is}.
\end{example}
If $g=nt(\sigma)$ is chosen uniformly at random, then $n$ and
$t(\sigma)$ are independent and uniform. Multiplying a sequence of
$g_i=n_it_i$ can be done by first multiplying the $t_i$, keeping
track of the carries, and then multiplying the $n_i$ and carries in
any order. Of course, here the carries are in $N$, not necessarily
binary.

Given $t_1,t_2,\dots,t_k$, we may form a two-column array with $t_i$ in
the first column, the successive remainders $r_1,r_2,\dots r_k$ in the
second column and ``carries'' $f_1,f_2,\dots f_{k-1}$ (elements in $N$)
placed in between.
\begin{lem}
If coset representatives $t_i,\ 1\leq i\leq k$, are chosen uniformly, then
\begin{enumerate}
\item the remainder process $r_i,\ 1\leq i\leq k$, is uniform and
  independent, and
\item the carries $f_i,\ 1\leq i\leq k-1$, form a stationary,
  one-dependent, two-block factor.
\end{enumerate}
\label{lem71}
\end{lem}
\begin{proof}\

\begin{enumerate}
\item Since $r_1=t_1$, $r_1$ is uniform. Successive $r_i$ are formed
  by multiplying $r_{i-1}$ by $t_i$. There is a unique choice of $t_i$
  giving $r_i$, so $r_i$ is uniform and independent of $r_1,\dots,
  r_{i-1}$.
\item Consider two successive remainders and the unique $t$ giving
  rise to them:
\begin{equation*}\begin{array}{cc}
&r_{i-1}\\
t_i&r_i\end{array}
\end{equation*}
Since $t_i$ is uniquely determined, $r_{i-1}t_i=r_if(r_{i-1},t_i)$ is
uniquely determined. It follows that the $f_i$ process is a two-block
process: generate $r_1,r_2\dots$ uniformly and independently, set
$f_i=h(r_i,r_{i+1})$, with $h(r,r')=f(r,t)$ where $r^{-1}r' \in tN$
determines $t$ uniquely. Because two-block processes are one-dependent,
this completes the proof.\qedhere
\end{enumerate}
\end{proof}
\begin{cor}
With the notation of the lemma, define a binary process
$B_1,\dots,B_{k-1}$ as
\begin{equation*}
B_i=\begin{cases}
1&\text{if }f_i\neq\emph{id}\\
0&\text{if }f_i=\emph{id}.
\end{cases}\end{equation*}
Then $\{B_i\}$ is a stationary, one-dependent, two-block factor.
\label{cor71}
\end{cor}

\begin{example} \label{ex72} As in Example \ref{ex71}, let $G=Q_8,\ N=\{\pm1\},
\ F=G/N\cong C_2\times C_2$. The natural choice of coset representatives
$\{1,i,j,k\}$ gives rise to the following two-block representation: choose
$U_1,U_2,\dots$ uniformly and independently in $\{1,i,j,k\}$. Let
$B_i=h(U_i,U_{i+1})$ with $h(1,x)=1$ (all $x$), $h(x,1)=-1 (x\neq1),\
h(i,j)=h(j,k)=h(k,i)=-1,\ h(i,k)=h(k,j)=h(j,i)=1,\ h(x,x)=1,\ (x\neq1)$.
\end{example}
\begin{example}(dihedral group)
Let $D_8=\langle x,y:x^2=y^2=(xy)^4=1\rangle$. If $z=xy$, this $8$-element
group has center $N=\{1,z^2=-1\}$. Choosing coset representatives
$1,x,y,z$, the cosets multiply as
\begin{equation*}\begin{array}{c|cccc}
&1&x&y&z\\\hline
1&1&x&y&z\\
x&x&1&z&y\\
y&y&-z&1&-x\\
z&z&-y&x&-1
\end{array}\end{equation*}
From this, elementary manipulations show that $P(i$ ones in a row)
$=1/4^i$. Thus, the carries process is independent with $P(B_i=1)=1/4, \
P(B_i=0)=3/4$ for all $i$.

We mention that $D_8$ can also be represented as the extension of the
normal subgroup $C_4$ by the factor group $C_2$. Here, coset
representatives for $C_4$ can be chosen so that there are no carries (i.e.
the extension ``splits''). \label{ex73}
\end{example}

\begin{example}(extensions of $C_2$ by $C_m$) A central extension of
$N=C_2$ by $C_m$ is abelian. It follows that when $m$ is odd, $G=C_{2m}$
is the only central extension of $N=C_2$ by $C_m$, and when $m$ is even,
there are two central extensions $C_2\times C_m$ and $C_{2m}$. The
extension $C_2\times C_m$ splits and choosing $\{(0,i), 0 \leq i \leq
m-1\}$ as the coset representatives, there are no carries. For $C_{2m}$,
with $C_2\cong\{0,m\}$, choose coset representatives $0,1,2,\dots,m-1$.
Thus
\begin{equation*}
f(i,j)=\begin{cases}
1&\text{if }i+j\geq m\\
0&\text{if }0\leq i+j<m.\end{cases}\end{equation*} As for usual addition,
there is a carry if and only if there is a descent in the remainder
column. Thus
\begin{equation*}
P(B_1=B_2=\dots=B_{i-1}=1)=\frac{\binom{m}{i}}{m^i}.
\end{equation*}
\label{ex74}
\end{example}

\begin{example} Let $G=C_2 \times C_2 \times C_2$ and let $N=\{(0,0,0),(1,1,1)\}$.
With coset representatives $(0,0,0),(1,0,0),(0,1,0),(1,1,0)$, there are
never carries, but with coset representatives
$(0,0,0),(1,0,0),(0,1,0),(0,0,1)$, there are carries. This example shows
that the one-dependent process $B_1,B_2,\cdots$ depends not only on $G$
and $N$, but also on the choice of coset representatives. \label{ex75}
\end{example}

We have not embarked on a systematic study of carries for finite groups
and believe that there is much more to do. We do note that by a result in
the next section (Theorem \ref{rdet}), the above processes, being
one-dependent, are determinantal.

The basic carries argument works for infinite groups as well. Let
$G$ be a locally compact group and $H$ a closed normal subgroup.
Suppose that $H$ is in the center of $G$ and that $G/H$ is compact.
Choose coset representatives $t(\sigma) \in G$ for $\sigma \in G/H$.
As in the finite case, these define factor sets $f(\sigma,\tau)$ by
$t(\sigma) t(\tau) = t(\sigma \tau) f(\sigma,\tau)$. Write $g=nt$.
Since $G/H$ is compact, it has an invariant probability measure.
Choosing $t_1,t_2,\cdots$ independently from this measure and
multiplying as above gives remainders and a carries process. Just as
above, the remainders are independent and uniformly distributed in
$G/H$ and the carries process (with values in $N$) is a
one-dependent, two-block factor process.

\begin{example} A lovely instance of this set-up explains a classical
identity. We begin with the motivation and then translate. Let $\sigma$ be
a permutation with number of descents $d(\sigma)$. It is known that, for
$\sigma$ chosen from the uniform distribution on $S_n$,
\begin{equation}
P\left( d(\sigma)=j \right) = P \left( j \leq U_1+ \cdots +U_n <j+1
\right). \label{918}
\end{equation} On the right, $U_1,U_2,\cdots,U_n$ are independent uniforms
on $[0,1]$. The density and distribution function for $U_1+\cdots+U_n$ was
derived by Laplace; see \cite{fel}. Foata \cite{fo} proved \eqref{918} by
combining Laplace's calculation with an identity of Worpitzky. Richard
Stanley \cite{sta2.5} gives a bijective proof involving an elegant
dissection of the $n$-dimensional hypercube. Jim Pitman \cite{Pi} gives
the following ``proof from the book'' which is a continuous version of our
carries argument from \ref{sec1}: form two columns

\begin{equation*}\begin{array}{rr}
U_1 & V_1\\
U_2 & V_2 \\
\cdot & \cdot \\
\cdot &\cdot \\
U_n & V_n
\end{array}\end{equation*} On the left are independent uniforms on $[0,1]$. On the
right are their remainders when added mod 1; so $V_1=U_1$, $V_2=U_1+U_2$
(mod $1$), $\cdots$. The $V_i$ are similarly independent uniforms on
$[0,1]$. Place a dot at position $i$ every time the partial sum
$U_1+\cdots+U_{i+1}$ crosses an integer. Call these dots carries. As in
the discrete case, there is a dot at position $i$ if and only if there is
a descent $V_{i+1}<V_i$. The number of dots is the integer part of
$U_1+\cdots+U_n$ and also the number of descents, proving \eqref{918}. Of
course, the distribution of the descent {\it process} is the same as the
carries process.

In the language of group theory, let $G=\mathbb{R}, N=\mathbb{Z}$, and
$G/N \cong S_1$, the circle group. Choose coset representatives as $[0,1)$
and factor sets in $\{0,1\}$.
\end{example}

\section{Proofs and generalizations} \label{sect7}

This section proves two of our main results for general
one-dependent processes (we do not assume stationarity in this
section). In Subsection \ref{all}, it is shown that all
one-dependent point processes on $\mathbb{Z}$ are determinantal, a
result which is new even in the stationary case. Subsection
\ref{class} proves that a point process $P$ on a finite set $\x$,
with $P$ given as a certain-shaped determinant, is one-dependent and
determinantal. This covers quite a few examples from previous
sections and is particularly useful in situations (such as Example
\ref{des2}) where the one-dependence is not apriori obvious.

\subsection{One dependent processes are determinantal} \label{all}

For a random point process on a discrete set $\x$, we define the
correlation function $\rho$ by \[ \rho(A) = P \{ S:S \supseteq A \}.\]
Then one dependence on (a segment of) $\mathbb{Z}$ is equivalent to the
condition that $\rho(X \cup Y) = \rho(X) \rho(Y)$ whenever $dist(X,Y) \geq
2$.

\begin{thm} \label{rdet} Any one-dependent point process on (a segment of
$\mathbb{Z}$) is determinantal. Its correlation kernel can be written in
the form $K(x,y)=$
\[  \begin{array}{ll}
0 & \mbox{if $x-y \geq 2$}\\
-1 & \mbox{if $x-y = 1$}\\
\sum_{r=1}^{y-x+1} (-1)^{r-1} \sum_{x=l_0<l_1<\cdots<l_r=y+1} \rho
\left([l_0,l_1) \right) \rho \left([l_1,l_2) \right) \cdots \rho
\left([l_{r-1},l_r) \right) & \mbox{if $x \leq y$}
\end{array} \] Here the notation $[a,b)$ stands for
$\{a,a+1,\cdots,b-1\}$. \end{thm}

For example, $K(x,x)=\rho(\{x\})$, $K(x,x+1)=\rho(\{x,x+1\})-\rho(\{x\})
\rho(\{x+1\})$, etc.

\begin{proof} By one-dependence, it is enough to verify that with $K$ as above, \[
\det[K(x+i,x+j)]_{i,j=0}^{y-x} = \rho \left([x,y+1) \right) \] for any $x
\leq y$.

We use induction on $y-x$. For $y=x$ the statement is trivial. Otherwise,
one has that \[ \det[K(x+i,x+j)]_{0}^{y-x} = \det \begin{bmatrix}
K(x,x)& K(x,x+1) & \dots & K(x,y)\\
-1     & K(x+1,x+1) &\dots &\vdots\\
      & -1     &\ddots&\vdots\\
      &      & -1    & K(y,y)
\end{bmatrix}. \]

When expanding this determinant, various numbers of $-1$'s from the
subdiagonal can be used. Observe that if we do not use the $-1$ in
position $(i,i-1)$, then we can compute the corresponding contribution,
because if we replace that $-1$ by $0$, the determinant splits into the
product of two each of which is computable by the induction hypothesis.

Similarly, if we insist on not using several $-1$'s, then the contribution
is the product of several determinants. Thus we obtain by
inclusion-exclusion that
\begin{equation} \begin{aligned}
\det[K(x+i,x+j)]_{0}^{y-x} & = K(x,y) + \sum_{x<l<y+1} \rho
\left([x,l)\right) \rho \left([l,y+1)\right) \\
& - \sum_{x<l_1<l_2<y+1} \rho \left([x,l_1)\right) \rho
\left([l_1,l_2)\right) \rho \left([l_2,y+1)\right) + \cdots.
\end{aligned} \label{eqv}
\end{equation} The $K(x,y)$ term corresponds to using all $-1$'s, the sum
over $l$ corresponds to not using a $-1$ at least at location
$(l,l-1)$, the sum over $l_1,l_2$ corresponds to not using a $-1$ at
least at locations $(l_j,l_j-1), j=1,2$, etc. Our definition of
$K(x,y)$ is such that the right hand side of \eqref{eqv} is equal to
$\rho \left([x,y+1)\right)$, as desired.
\end{proof}

To state some corollaries of Theorem \ref{rdet}, we use the concept of
particle-hole involution. Essentially, given a subset $\cal{N}$ of $\x$,
the involution maps a point configuration $S \subset \x$ to $S
\bigtriangleup \cal{N}$ (here $\bigtriangleup$ is the symbol for symmetric
difference). This map leaves intact the particles of $S$ outside of
$\cal{N}$, and inside $\cal{N}$ it loses the particles and picks up the
``holes'' (points of $\cal{N}$ free of particles).

\begin{cor} The class of one-dependent processes is closed under the
operations of particle-hole involution on any fixed subset, intersections
of independent processes, and unions of independent processes. \end{cor}

\begin{proof} For intersections, the claim follows from the definitions.

For particle-hole involutions, note that the property of being
one-dependent follows from the kernel begin 0 on the second subdiagonal
and below. On the other hand, for determinantal point processes the
``complementation principle'' \citep[Sec.~A.3]{boo} says that the
particle-hole involution can be implemented by the following change in the
kernel:
\[ \begin{bmatrix}
A& B \\
C& D
\end{bmatrix} \rightarrow \begin{bmatrix}
A& B \\
-C& I-D
\end{bmatrix}, \] where the block structure corresponds to the splitting
into the noninverted and inverted parts. Clearly, this keeps the property
of having $0$'s on the second subdiagonal and below intact.

Finally, unions can be reduced to intersections by the particle-hole
involution. \end{proof}

{\it Remark}. The class of determinantal point processes is {\it not}
closed under intesections/unions.

\vspace{.1in}

Let us now look at the translation invariant case. Then $\rho \left([x,y)
\right) = \rho_{y-x}$ is a function of $y-x$ only, and $\rho_k$ is the
chance of $k$ consecutive ones. Set
\[ R(z) = 1+z+\sum_{k \geq 1} \rho_k z^{k+1}.\]

\begin{cor} \label{cor73} In the translation invariant case, the kernel $K(x,y)=k(y-x)$
is also translation invariant, and \[ \sum_{n \in \mathbb{Z}} k(n)z^n =
\frac{1}{1-R(z)}.\] \end{cor}

\begin{proof}
\begin{eqnarray*}
\frac{1}{1-R(z)} & = & - \frac{1}{z} \cdot \frac{1}{1+\rho_1 z+\rho_2 z^2
+
\cdots} \\
& = & -\frac{1}{z} \left[1 - \sum_{m \geq 1} \rho_m z^m + \left(\sum_{m
\geq 1} \rho_m z_m \right)^2 - \cdots \right] \\
& = & - \frac{1}{z} + \frac{1}{z} \left[ \sum_{m \geq 1} \rho_m z^m -
\left(\sum_{m \geq 1} \rho_m z^m \right)^2 + \cdots \right],
\end{eqnarray*} which agrees with the formula of Theorem \ref{rdet}.
\end{proof}

The generating function $R(z)$ also behaves well with respect to the
particle-hole involution on $\mathbb{Z}$.

\begin{prop} The particle-hole involution on $\mathbb{Z}$ with
$\cal{N}=\mathbb{Z}$ replaces $R(z)$ by $1/R(-z)$. \end{prop}

\begin{proof} We have \[ \frac{1}{1-1/R(-z)} = \frac{R(-z)}{R(-z)-1} = 1 -
\frac{1}{1-R(-z)}.\] Hence by Corollary \ref{cor73}, changing $R(z)$ to
$\tilde{R}(z)=1/R(-z)$ leads to the following change in the correlation
kernel: \[ \tilde{k}(n)= \delta_{0,n} - (-1)^n k(n).\] This is equivalent
to $\tilde{K}(x,y) = \delta_{x,y}- \frac{(-1)^x}{(-1)^y}K(x,y)$. This is
the same as $K \rightarrow I-K$, which corresponds to the particle-hole
involution \cite[Sec.~A.3]{boo}. \end{proof}

\begin{example} The descent process on $\mathbb{Z}$ (Example \ref{des1})
corresponds to $\rho_m=\frac{1}{(m+1)!}$; thus $R(z)=e^z$. The
particle-hole involution is given by $\tilde{R}(z)=1/e^{-z}=e^z$, which is
the same.
\end{example}

\begin{example} \label{bor2} The intersection of $r$ independent descent
processes on $\mathbb{Z}$ corresponds to $\rho_m= \frac{1}{(m+1)!^r}$,
hence \[ R_{\cap}^{(r)}(z) = \sum_{m \geq 0} \frac{z^m}{(m!)^r}.\] The
correlation kernel is given by $\sum_{n \in \mathbb{Z}} k(n)z^n =
\frac{1}{1-R_{\cap}^{(r)}(z)}$.
\end{example}

\begin{example} The union of $r$ independent descent processes is the
particle-hole involution of Example \ref{bor2}. Thus
$R_{\cup}^{(r)}(z)=1/R_{\cap}^{(r)}(-z)$ and \[ \sum_{n \in \mathbb{Z}}
k(n) z^n = \frac{1}{1-1/R_{\cap}^{(r)}(-z)}.\] Replacing $z$ by $-z$
doesn't affect correlations, and we recover Example \ref{int}.
\end{example}

\subsection{A class of determinantal processes} \label{class}

For any $n=2,3,\dots$, consider a probability  measure $P_n$ on all
subsets $S=\{s_1<s_2<\dots<s_k\}\subseteq[n-1]$ given by
\begin{equation}
P_n(S)=h(n)\,\text{det}\left[e(s_i,s_{j+1})\right]_{i,j=0}^k, \label{51}
\end{equation}
for some $h:\mathbb{N}\to\mathbb{C}$ and
$e:\mathbb{N}\times\mathbb{N}\to\mathbb{C}$ with the notation $s_0=0,
s_{k+1}=n$. We assume that $e(i,j)=0$ for $i>j$, and that $e(i,i)=1,
e(i,i+1)>0$ for all $i$.
\begin{thm}
If \eqref{51} holds for some fixed $n$ with $e(i,j)=0$ for $i>j$ and
$e(i,i)=1,e(i,i+1)>0$ for all $i$, then $P_n$ is a determinantal,
one-dependent process with correlation functions
\begin{equation*}
\rho(A)=P_n\{S:S\supseteq A\}=\emph{det}\left[K(a_i,a_j)\right]_{i,j=1}^m
\qquad\text{for }A=\{a_1,a_2,\dots,a_m\},
\end{equation*}
with correlation kernel
\begin{equation} \label{kform}
K(x,y)=\delta_{x,y}+(E^{-1})_{x,y+1},
\end{equation}
where $E$ is the upper triangular matrix $E=[e(i-1,j)]_{i,j=1}^n$,
\begin{equation*}
E=\begin{bmatrix}
e(0,1)&e(0,2)&\dots &e(0,n)\\
      &e(1,2)&\dots &e(1,n)\\
      &      &\ddots&\vdots\\
      &      &      &e(n-1,n)
\end{bmatrix}\ .
\end{equation*}
In addition, $h(n)=(\det\,E)^{-1}=(e(0,1)e(1,2)\dots e(n-1,n))^{-1}$.
\label{thm51}
\end{thm}
\begin{cor}
Assume further that $e(i,j)=e(j-i)$. Then the point process is stationary.
If $\hat{e}(z)=\sum_{l=0}^\infty e(l)z^l$, then $K(x,y)=k(y-x)$ and
\begin{equation*}
\hat{k}(z):=\sum_{l=-\infty}^\infty k(l)z^l=\dfrac1{1-1/\hat{e}(z)}.
\end{equation*}
\label{cor1}
\end{cor}
\begin{proof} [Proof of Theorem \ref{thm51}] Set $L=[e(i-1,j)+\delta_{i-1,j}]_{i,j=1}^n$.
Thus $L$ appears as
\begin{equation*}
L=\begin{bmatrix}
e(0,1)&e(0,2)&\dots &e(0,n)\\
1     &e(1,2)&\dots &\vdots\\
      &1     &\ddots&\vdots\\
      &      &1     &e(n-1,n)
\end{bmatrix}.
\end{equation*}
For any function $f:\{1,\dots,n-1\}\to\mathbb{C}$, by \eqref{51},
\begin{align*}
& E\left(\prod_{s_i \in S} f(s_i)\right)\\ &=
\sum_{\substack{0<s_1<\dots<s_k<n\\k=0,\dots,n-1}}P_n\left(\{s_1<s_2<\dots<s_k\}
\right)f(s_1)\dots f(s_k)\\
&=h(n)\sum_{\substack{0<s_1<\dots<s_k<n\\k=0,\dots,n-1}}\text{det}
\left[L\binom{1,s_1+1,\dots,s_k+1}{s_1,\dots,s_k,n}\right]f(s_1)\dots f(s_k)\\
&=h(n) \cdot \text{det} \begin{bmatrix}
f(1)e(0,1)& f(2) e(0,2)& \dots & f(n-1) e(0,n-1) & e(0,n)\\
f(1)-1  & f(2) e(1,2)& \dots & f(n-1) e(1,n-1) & e(1,n)\\
      &    f(2) - 1  & \ddots & f(n-1) e(2,n-1) & e(2,n)\\
      &      &    \ddots   & \dots & \dots \\
      & & & f(n-1) -1 & e(n-1,n)
\end{bmatrix}\\
&=h(n) \cdot\text{det}\left[\begin{bmatrix}
0&&&\\
-1&\ddots&&\\
&\ddots&\ddots&\\
&&-1&0
\end{bmatrix}+L\begin{bmatrix}
f(1)&&&\\
&\ddots&&\\
&&f(n-1)&\\
&&&1
\end{bmatrix}\right].
\end{align*} This is the generating functional of $P_n$. In the second equality,
the determinant is of the minor of $L$ with rows $1,s_1+1,\cdots,s_k+1$
and columns $s_1,\cdots,s_k,n$. In the third inequality, the $2^{n-1}$
possible summands correspond to choosing which of the first $n-1$ matrix
columns use the $-1$ coming from $f(i)-1$ in the determinant expansion.

If $f=1+g$, the generating functional can be expressed in terms of the
correlation functions:
\begin{eqnarray*}
& &
1+\sum_{\substack{s_1<\dots<s_m\\m=1,2,\dots}}\rho_m(s_1,\dots,s_m)g(s_1)\dots
g(s_m)\\ &
=  & E\left[\prod_{s_i \in S} \left(1+g(s_i)\right)\right]\\
&= & h(n)\,\text{det}\left(E+L\begin{bmatrix}
g(1)&&&\\
&\ddots&&\\
&&g(n-1)&\\
&&&0
\end{bmatrix}\right)\\
&= & h(n)\,\text{det}(E)\cdot\text{det}\left(I+E^{-1}L\cdot\begin{bmatrix}
g(1)&&&\\
&\ddots&&\\
&&g(n-1)&\\
&&&0
\end{bmatrix}\right).
\end{eqnarray*}
This holds for all $g$. First take $g=0$ to see $h(n)\,\text{det}(E)=1$.
Next note the expansion: if $M$ is $n\times n$, then
$\text{det}[I+M]=\sum_S \text{det}(M(S))$ with the sum over all $2^n$
subsets of $[n]$, and $M(S)$ the minor with rows and columns in $S$. This
gives that $\rho_m(s_1,\dots,s_m)=\text{det}[K(s_i,s_j)]_{i,j=1}^m$, where
$K(x,y)=(E^{-1}L)_{xy}$. Then
\begin{equation*}
E^{-1}L=E^{-1}\left(\begin{bmatrix}
0&&&\\
1&\ddots&&\\
&\ddots&\ddots&\\
&&1&0
\end{bmatrix}+E\right)=I +E^{-1} \begin{bmatrix}
0&&&\\
1&\ddots&&\\
&\ddots&\ddots&\\
&&1&0
\end{bmatrix},\qedhere
\end{equation*} and the proof of \eqref{kform} is complete.

Finally, we note that since $K(x,y)$ vanishes below the first subdiagonal,
$P_n$ is one-dependent. \end{proof}
\begin{proof}[Proof of Corollary]
Here $E$ is a Toeplitz matrix with symbol $(\hat{e}(z)-1)/z$. It is
triangular, so $E^{-1}$ is a Toeplitz matrix with symbol
$z/(\hat{e}(z)-1)$. Thus $K$ is a Toeplitz matrix with symbol
$\hat{k}(z)=1+\tfrac1{\hat{e}(z)-1}=\tfrac1{1-(\hat{e}(z))^{-1}}$.
\end{proof}
We also note that in the translation invariant case, there is an
expression for $P_n(S)$ in terms of skew Schur functions of ribbon type.

\begin{cor} \label{skew2} As in Corollary \ref{cor1}, assume further that
$e(i,j)=e(j-i)$. Let $\lambda$ and $\mu$ be the partitions defined by
\[ \lambda_i=n-s_{i-1}-k+i-1 \ , \ \mu_i=n-s_{i}-k+i-1, \ \ 1 \leq i
\leq k+1.\] Let $\lambda',\mu'$ denote the transpose partitions of
$\lambda$ and $\mu$ and let $s_{\lambda'/\mu'}$ denote the corresponding
skew-Schur function, obtained by specializing the elementary symmetric
functions $e_i$ to equal $e(i)$. Then $P_n(S) = \frac{1}{e(1)^n}
s_{\lambda'/\tau'}$. \end{cor}

\begin{proof} By assumption, formula \eqref{51} becomes
\[ P_n(S)=h(n) \cdot \text{det}\left(e_{s_{j+1}-s_i}\right)_{i,j=0}^k.\]
 Now use the argument of Theorem \ref{skewschur}, together with the
identification of the normalizing constant $h(n)=1/e(1)^n$ in Theorem
\ref{thm51}. \end{proof}

\begin{rems}\

\begin{enumerate}

\item Another approach to Theorem \ref{thm51} is via the theory of conditional
$L$-ensembles in Borodin--Rains \cite[Prop.~1.2]{br}.

\item In the translation invariant case $e(i,j)=e(j-i)$, Stanley
\cite[p.~90, Ex.~14]{sta2} shows that
\begin{equation*}
P_n\left([n-1]\right)=\text{chance that all sites are occupied}
\end{equation*}
is the coefficient of $z^{n}$ in the power series $h(n)\hat{e}(-z)^{-1}$.
To prove this using symmetric function theory, note from Corollary
\ref{skew2} that
\[
P_n\left([n-1]\right)=h(n) \cdot s_{(n)}.\] Here the Schur function
$s_{(n)}=h_n$ (where $h_n$ is the $n$th complete homogeneous symmetric
function and $h(n)$ is the normalizing constant in \eqref{51}). If
$E(z)=\sum_{r=0}^\infty e_r z^r$ is the generating function for elementary
symmetric functions and $H(z)=\sum_{r=0}^\infty h_{r}z^r$, Macdonald
\cite[p.~21, (2.6)]{mac} shows that $H(z)E(-z)=1$. This gives Stanley's
formula since in Corollary \ref{skew2}, the value of $e_r$ is $e(r)$.
\end{enumerate}
\end{rems}
\begin{example} \textbf{(Descents in a random sequence)}\quad
From Fact 5 of \ref{sec2}, we are in the Toeplitz case with
$e(j)=\binom{j+b-1}{b-1}, \hat{e}(z)=(1-z)^{-b}$. Applying Corollary
\ref{cor1} yields Theorem \ref{thm31}. \label{exam51}
\end{example}
\begin{example} \textbf{(Descents in a uniform permutation)}\quad
From MacMahon's formula \eqref{41}, we are again in the Toeplitz case with
$e(j)=1/j!, \hat{e}(z)=e^{z}$. As in Example \ref{ex31}, one can replace
$z$ by $-z$ without changing determinants or the correlations functions.
This proves Theorem \ref{thm41}. \label{exam52}
\end{example}
\begin{example} \textbf{(Descents in a non-uniform permutation)}\quad
From Stanley's formula (part a of Proposition \ref{prop1}), we are again
in the Toeplitz case with $e(j)=1/j!_q$ ($q$-factorial). An identity of
Euler allows one to write
\[ \hat{e}(z) = \prod_{m \geq 0} \frac{1}{1-z(1-q)q^m} \] when $0<q<1,|z|<1$,
but the elementary description of the correlation functions given in
Proposition \ref{prop1} seems more useful. \label{exam53}
\end{example}

\section*{Acknowledgments} Borodin was partially supported by NSF grant
DMS 0707163. Diaconis was partially supported by NSF grant DMS 0505673.
Fulman was partially supported by NSF grant DMS 0802082 and NSA grant
H98230-08-1-0133.

\end{document}